\documentclass[runningheads]{llncs}
\usepackage{graphicx}
\usepackage[super]{nth}
\usepackage{subfigure}
\usepackage{amsmath}
\usepackage{algorithm}
\usepackage{algpseudocode}
\begin{document}
\title{Tabu Search and Simulated Annealing
metaheuristic algorithms applied to the RoRo
vessel stowage problem}
%
%
\titlerunning{Survey on Stowage Planning}

\author{Eghbal Hosseini}
\authorrunning{ Eghbal Hosseini}
%

\institute{Department of People and Technology Roskilde University, Denmark \\
\email{kseghbalhosseini@gmail.com, hosseini@ruc.dk}\\
}
\maketitle              
\begin{abstract}
The search heuristics Tabu search and Simulated annealing
are commonly used meta-heuristics. The two heuristics have different
ways of ensuring diversification. The heuristics can be implemented for
solving the stowage planning problem. The stowage planning problem
occurs every time a vessel is loaded. The idea is to stow the cargo in an
optimal manner satisfying a set of constraints and specifically the stability constraints. An optimal plan can be to assign the cargo to spots on
the vessel so that the vessel trim and draft are optimize ensuring a low
fuel consumption. Here the problem considered is the stowage planning
of trailers on a Roll-on/Roll-off vessel. Roll-on/Roll-off vessels carries
vehicles or trailers on the different levels of the vessel. When making a
stowage plan of a vessel ballast tanks can be adjusted to improve stability and to change draft. This leads to stability constraints which can
be very complex to model and solve using a mixed integer model and
solver. Complicating constraints occur in many different applications and
different forms and are the cause of the popularity of search algorithms,
such as tabu search and simulated annealing, for solving real-life applications. Results of the two meta-heuristics are shown for real-life stowage
planning cases for Roll-on/Roll-off vessels.

\end{abstract}

\section{Introduction}

\section{Modelling}
The problem is stowage planning for a ship with some unusable spaces. Cargos should be loaded in defined ports and unloaded in others. The objective of the problem is minimizing total time of load and unload of cargos. The problem has been modeled as a mixed-integer linear programming (MILP) problem as follows:

\begin{equation}
\label{32dddd1987}
\begin{aligned}
&\min && \sum_{i=1}^{k} \sum_{j=1}^{m\times n} t_{ij} 
\end{aligned}
\end{equation}
\begin{equation}
\label{32dddd19287}
\begin{aligned}
\sum_{i=1}^{k} x_{ij} \leq1,\\
\end{aligned}
\end{equation}
\begin{equation}
\begin{aligned}
\sum_{j=1}^{m\times n} x_{ij} =1, &\text{for i=1,2,...,k},\\
\end{aligned}
\label{100}
\end{equation}
\begin{equation}
\begin{aligned}
\sum_{i=1}^{k} y_{ip} \cdot  w_{i}\leq  w_{p}, &\text{for p=1,2,...,l},\\
\end{aligned}
\label{111}
\end{equation}
\begin{equation}
\begin{aligned}
t_{ij} =\left\{
        \begin{array}{c c}
            T &\text{i,j belong to a same category}\\
            T+Q &\text{Otherwise}
        \end{array}\right.
\end{aligned}
\label{111}
\end{equation}
\begin{equation}
\begin{aligned}
x_{ij} =\left\{
        \begin{array}{c c}
            1 &\text{Cargo i set in the cell j}\\
            0 &\text{Otherwise}
        \end{array}\right.
\end{aligned}
\label{111}
\end{equation}
\begin{equation}
\begin{aligned}
y_{ip} =\left\{
        \begin{array}{c c}
            1 &\text{Cargo i set in the deck p}\\
            0 &\text{Otherwise}
        \end{array}\right.
\end{aligned}
\label{111}
\end{equation}

The objective of the model, (1), minimises the total time of unloading cargos. Constraints (2) ensure that each cargo is loaded exactly once at a cell. Constraints (3) make sure that each cell will only have at most one cargo loaded. Constraints (4) make sure that the total weight of cargoes loaded on each deck does not exceed the maximum weight limit per deck. $w_{i}$ is the weight of the ith cargo. 
\section{LCA, VEA, SA, and TS Algorithms}
\subsection{Simulating Annealing (SA)}
There are several heuristic and meta-heuristic algorithms in the literature [1-23]. Simulated annealing (SA) is a meta-heuristic technique for approximating the global optimum of a given optimization problem. Specifically, it is a probabilistic method to approximate global optimization in a big search region for an optimization problem. It is often used for discrete optimization problem. 
In the Simulated Annealing case, the equation has been altered to the following:

\begin{equation}
\begin{aligned}
p =\left\{
        \begin{array}{c c}
            1 &{\Delta c \leq 0}\\
            e^{\frac{-\Delta c}{t} } &{\Delta c > 0}
        \end{array}\right.
\end{aligned}
\label{111}
\end{equation}

Where the delta c is the change in cost and the t is the current temperature. The P calculated in this case is the probability that we should accept the new solution. In our implementation the formulation of SA has been little changed based on the problem:

\begin{equation}
\begin{aligned}
p =\left\{
        \begin{array}{c c}
            e^{\frac{1}{\Delta c} } &{\Delta c \leq 0}\\
            \frac{e^{\frac{-\Delta c}{t} }}{M} &{\Delta c > 0}
        \end{array}\right.
\end{aligned}
\label{111}
\end{equation}

Which $M>1$. The process of the modified SA are as follows:\\
1.	The algorithm starts with an Initial Solution (IS).\\
2.	The population is generated close to initial solution.\\
3.	Each solution of population gets a probability based on (2).\\
4.	The solution with the biggest probability will be selected as the best solution of current population.\\
5.	 Set IS = best solution and go back to step 2. \\

\begin{algorithm}
\caption{ SA Procedure for Stowage Planning Problem} 
\label{kadds321aaafif}
\begin{algorithmic}[1]
\item n: Number of solutions
\item N: Number of Iterations
\item $\Delta c_{k}$: Difference of objective functions between solution k and initial solution
\item $\alpha$: A given positive number (less than size of the problem)
\item Generate a random initial feasible solution X0
\item Generate initial population near initial solution
\For{$i \gets 1$ to $N$}  
\For{$k \gets 1$ to $n$}  
\If{Xk is better than X0} 
\item $P_{k}=e^{\frac{1}{\Delta ck} } $
\ElsIf
\item $P_{k}=\frac{e^{\frac{-\Delta c}{t} }}{M} $
\EndIf
\EndFor
\item Find Maximum of Pk 
\item $P_{max}=Maximum of Pk$ and $r=k$
\item $Xbest=Xr$
\EndFor
\item $X0=Xbest$ and back to 6
\end{algorithmic}
\end{algorithm}

\subsection{Tabu Search (TS)}
Tabu Search (TS) is a metaheuristic search method employing local search methods used for mathematical optimization. Local searches take a feasible solution to the optimization problem and check its very close neighbors to find a better solution. Local search methods tend to become stuck in local optimal solutions. TS modifies the performance of local search by relaxing its basic rule. At each step of TS algorithm, changing of moves can be accepted if better solution is not available such as when the search is stuck at a strict local optimal. In addition, prohibitions (tabu) are introduced to discourage the search from coming back to previous visited solutions.

\begin{algorithm}
\caption{ TS Procedure for Stowage Planning Problem} 
\label{kadds321aaafif}
\begin{algorithmic}[1]
\item n: Number of solutions
\item s: Size of the problem
\item Rand: Random integer number between 1 and s
\item $\lambda$: A given integer positive number 
\item Generate initial population 
\For{$t \gets 1$ to $\lambda$}  
\For{$k \gets 1$ to $n$}
\item $Xk=Xk+Rand *\frac{Xk}{||Xk||}$ (Distribution of solutions in feasible region)
\EndFor
\item Find best solution 
\item Let $Xt= xbest$
\EndFor
\For{$t \gets 1$ to $\lambda$}  
\For{$i \gets 1$ to $n$}
\item $Xi=Xt+Rand *\frac{Xk}{||Xk||}$ (Distribution of best solutions)
\EndFor
\item Find best solution 
\item Let $Yt= xbest$
\EndFor
\end{algorithmic}
\end{algorithm}

\subsection{Laying Chicken Algorithm (LCA)}
Laying Chicken Algorithm has been inspired from behavior of laying hens. It focuses on finding an answer to how does the hen convert the egg to the chicken? LCA converts the feasible solutions to the optimal solution, same as what laying hen does from the eggs to the chickens. In fact, each feasible solution of a continuous programming problem displays an egg and the optimal solution of the problem is a chicken. Hens try to warm their eggs; this concept has inspired by LCA to change and improve the feasible solutions. 
There are the following steps to formulate of the behavior the hen in the  LCA optimizer [4]:
\begin{enumerate}
    \item The first egg which displays initial solution.
    \item More eggs displays initial population close to the initial solution.
    \item Improve solutions of population inspiring from warming eggs.
    \item Little mutation of solutions inspiring of rotation of eggs.
\end{enumerate}

\begin{algorithm}
\caption{ LCA Procedure for Stowage Planning Problem} 
\label{kadds321aaafif}
\begin{algorithmic}[1]
\item n: Number of solutions
\item N: Number of Iterations
\item $\alpha$: A given positive number (less than size of the problem)
\item Generate a random initial feasible solution X0
\item Generate initial population near initial solution
\For{$i \gets 1$ to $N$}  
\For{$k \gets 1$ to $n$}  
\If{Xk is not better than X0} 
\item $Xk=X0+\alpha *(\frac{Xk}{||Xk||})$
\EndIf
\item $Xk=Xk+\frac{Xk}{||Xk||}$
\EndFor
\EndFor
\end{algorithmic}
\end{algorithm}

\subsection{Volcano Eruption Algorithm (VEA)}
The Volcano Eruption Algorithm has been inspired from the nature of a volcano eruption. VEA optimizer imitates the process of volcano eruption, which is a hole on the earth’s surface. This phenomenon acts as a vent for release of pressurized gases, molten rock or magma deep beneath the surface of earth. Magma is passed through a channel from deep underground called the volcanic pipe. Magma erupts out of the earth’s surface when it reaches the hole on the surface. 
There are the following steps leading to formation of a volcano to VEA optimizer [5]:
\begin{enumerate}
    \item Rise of magma through the volcanic pipe.
    \item Volcanic eruption by rising of magma to the surface of the earth.
    \item Lava’s cooling down and therefore formation of a crust.
    \item Repetition of this process over time leading to several layers of rock that builds up over time resulting in a volcano.
\end{enumerate}

\begin{algorithm}
\caption{ VEA Procedure for Stowage Planning Problem} 
\label{kadds321aaafif}
\begin{algorithmic}[1]
\item n: Number of solutions
\item s: Size of the problem
\item Rand: Random integer number between 1 and s
\item $\lambda$: A given integer positive number 
\item Generate initial population 
\For{$t \gets 1$ to $\lambda$}  
\For{$k \gets 1$ to $n$}
\item $Xk=Xk+Rand *\frac{Xk}{||Xk||}$ 
\EndFor
\item Find best solution 
\item Let $Xt= xbest$
\EndFor
\For{$t \gets 1$ to $\lambda$}  
\For{$i \gets 1$ to $n$}
\item $Xi=Xt+Rand *\frac{Xk}{||Xk||}$ 
\EndFor
\item Find best solution 
\item Let $Yt= xbest$
\EndFor
\end{algorithmic}
\end{algorithm}
\subsection{Multiverse Algorithm (MVA)}
\begin{algorithm}
\caption{ MVA Procedure for Stowage Planning Problem} 
\label{kadds321aaafif}
\begin{algorithmic}[1]
\item N: Number of solutions
\item s: Size of the problem
\item Rand: Random integer number between 1 and s
\item $\lambda$: A given integer positive number 
\item m: Number of solutions for the next population 
\item Generate initial population
\For{$t \gets 1$ to $\lambda$}  
\For{$l \gets 1$ to $N$}   
\For{$k \gets 1$ to $m$}
\item $Xk=Xk+Rand *\frac{Xk}{||Xk||}$ (Distribution of solutions in feasible region)
\EndFor
\item Find best solution 
\item Let $Xt= xbest$
\EndFor
\EndFor
\For{$t \gets 1$ to $\lambda$}  
\For{$i \gets 1$ to $m$}
\item $Xi=Xt+Rand *\frac{Xk}{||Xk||}$ (Distribution of best solutions)
\EndFor
\item Find best solution 
\item Let $yt= xbest$
\EndFor
\end{algorithmic}
\end{algorithm}

Algorithms 1-5 show the pseudocodes of SA, TS, LCA, VEA, and MVA algorithms respectively.\\

\section{Computational Results}
In this section, four algorithms: Tabu Search (TS), Simulated Annealing (SA), Laying Chicken Algorithm (LCA) and Volcano Eruption Algorithm (VEA) are implemented to solve some benchmarks of the problem. 

\begin{figure}
    \centering
    \subfigure[Unusable space (dark blue) before loading]{\includegraphics[width=0.49\textwidth]{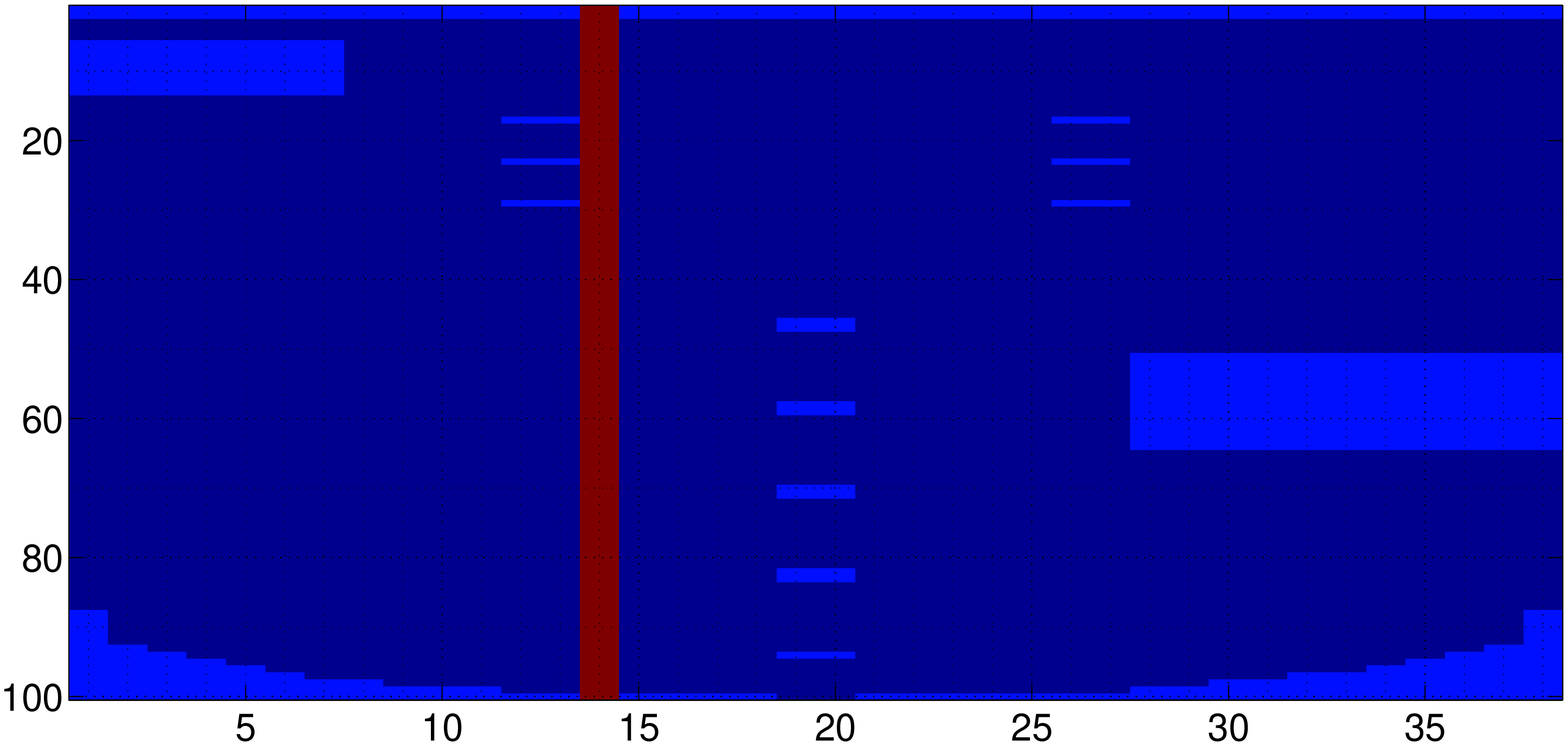}} 
    \subfigure[Initial solution by LCA]{\includegraphics[width=0.49\textwidth]{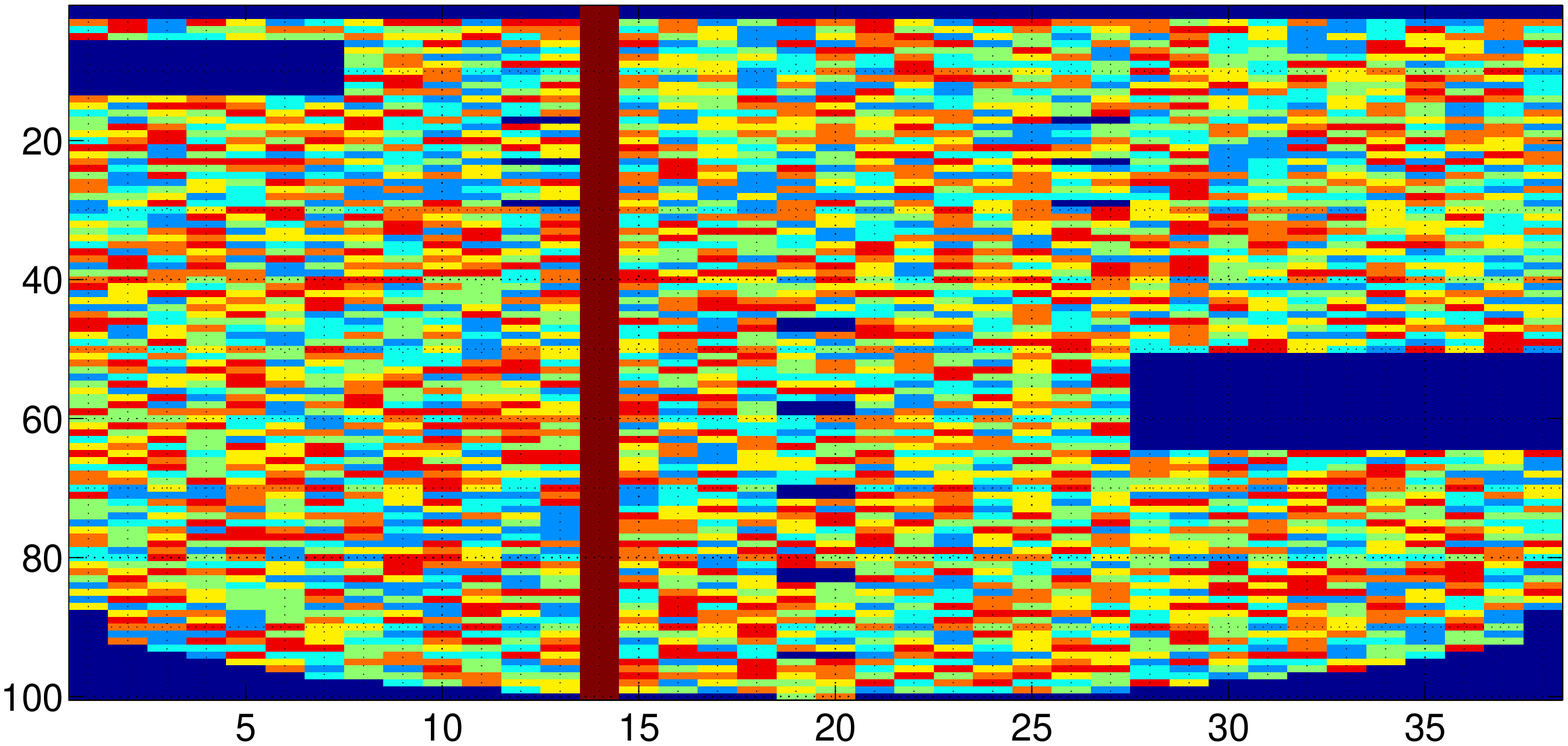}}
    \caption{Matrix format of ship before loading and initial solution by LCA}
    \label{fig:foobar}
\end{figure}

\begin{figure}
    \centering
    \subfigure[Problem 1]{\includegraphics[width=0.49\textwidth]{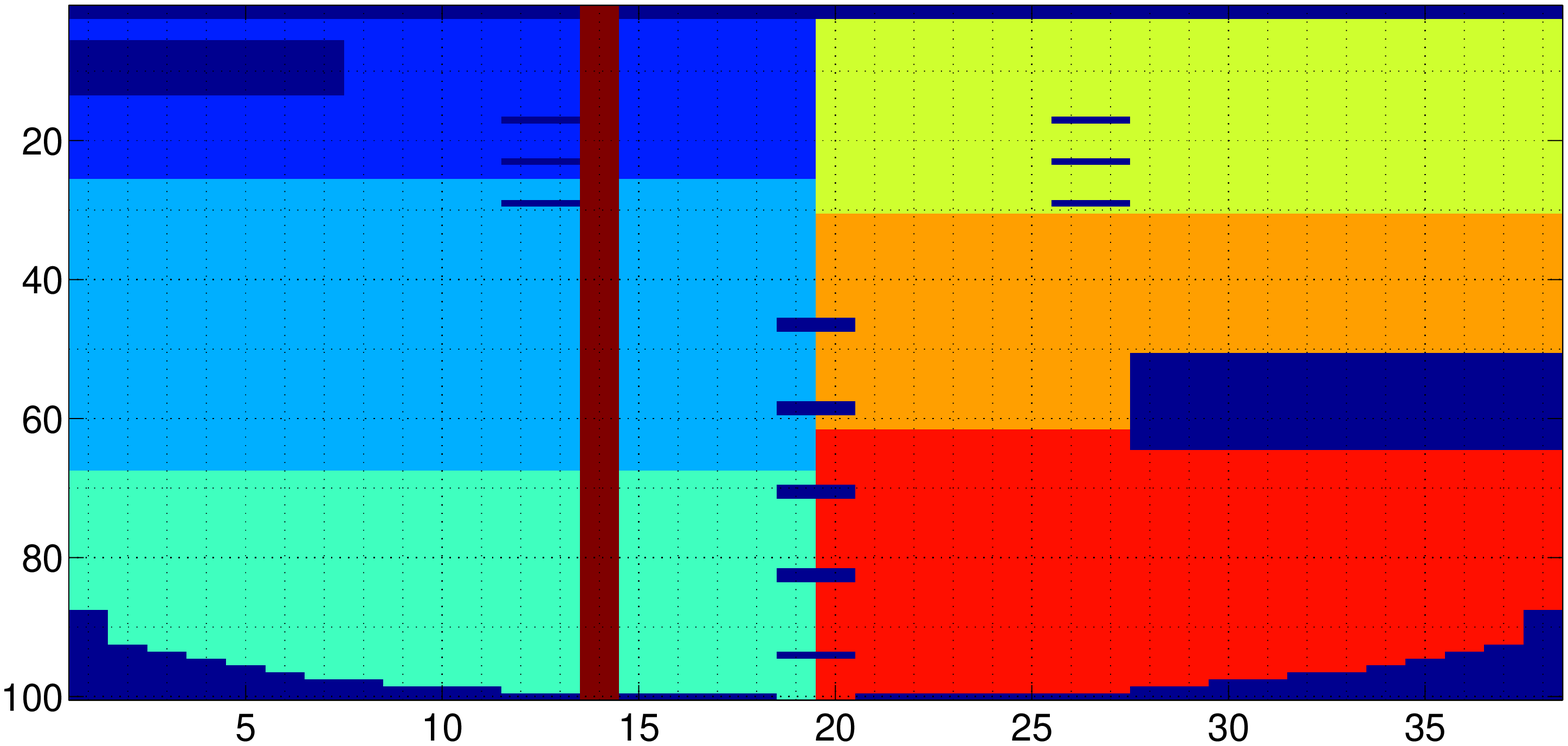}} 
    \subfigure[Problem 2]{\includegraphics[width=0.49\textwidth]{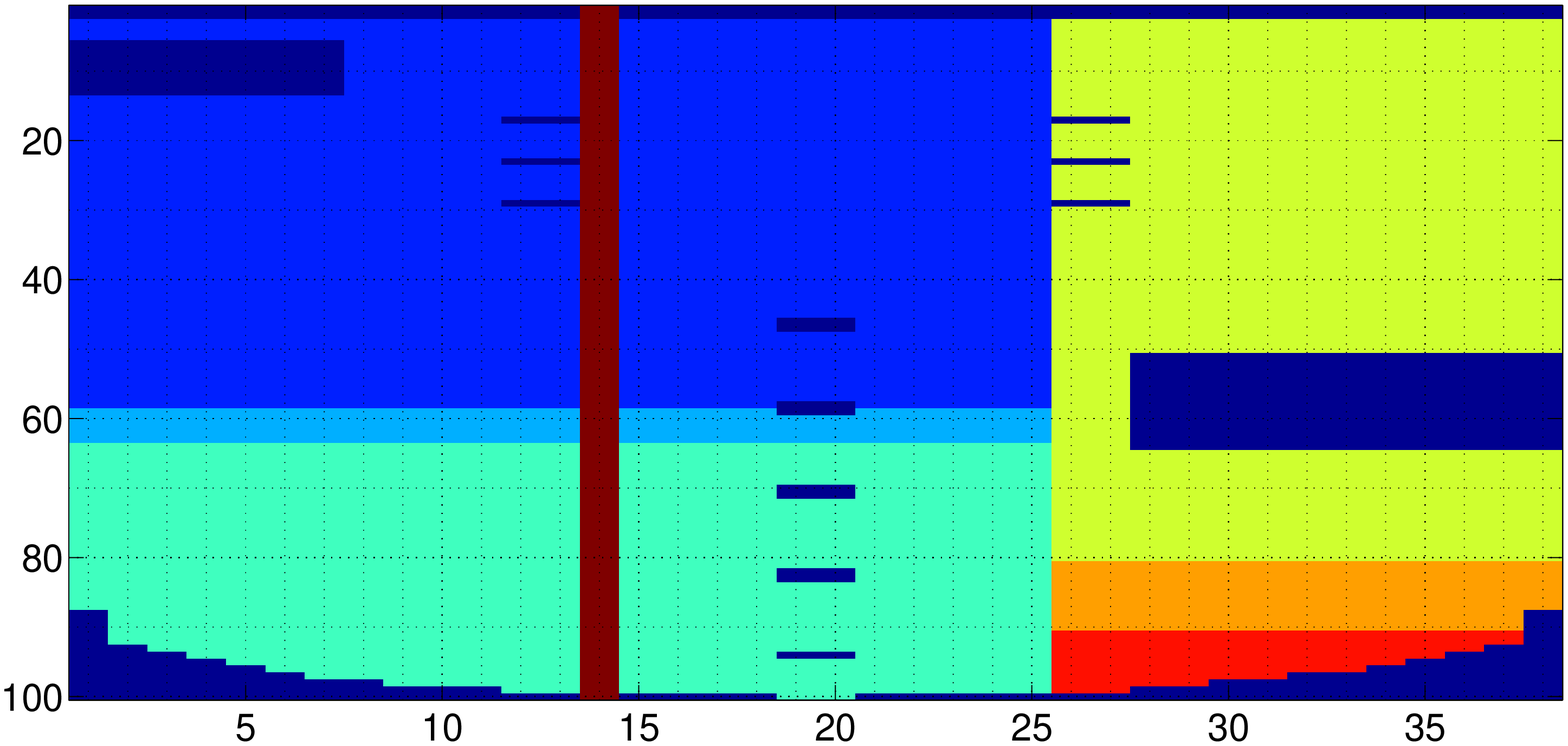}}
    \subfigure[Problem 3]{\includegraphics[width=0.49\textwidth]{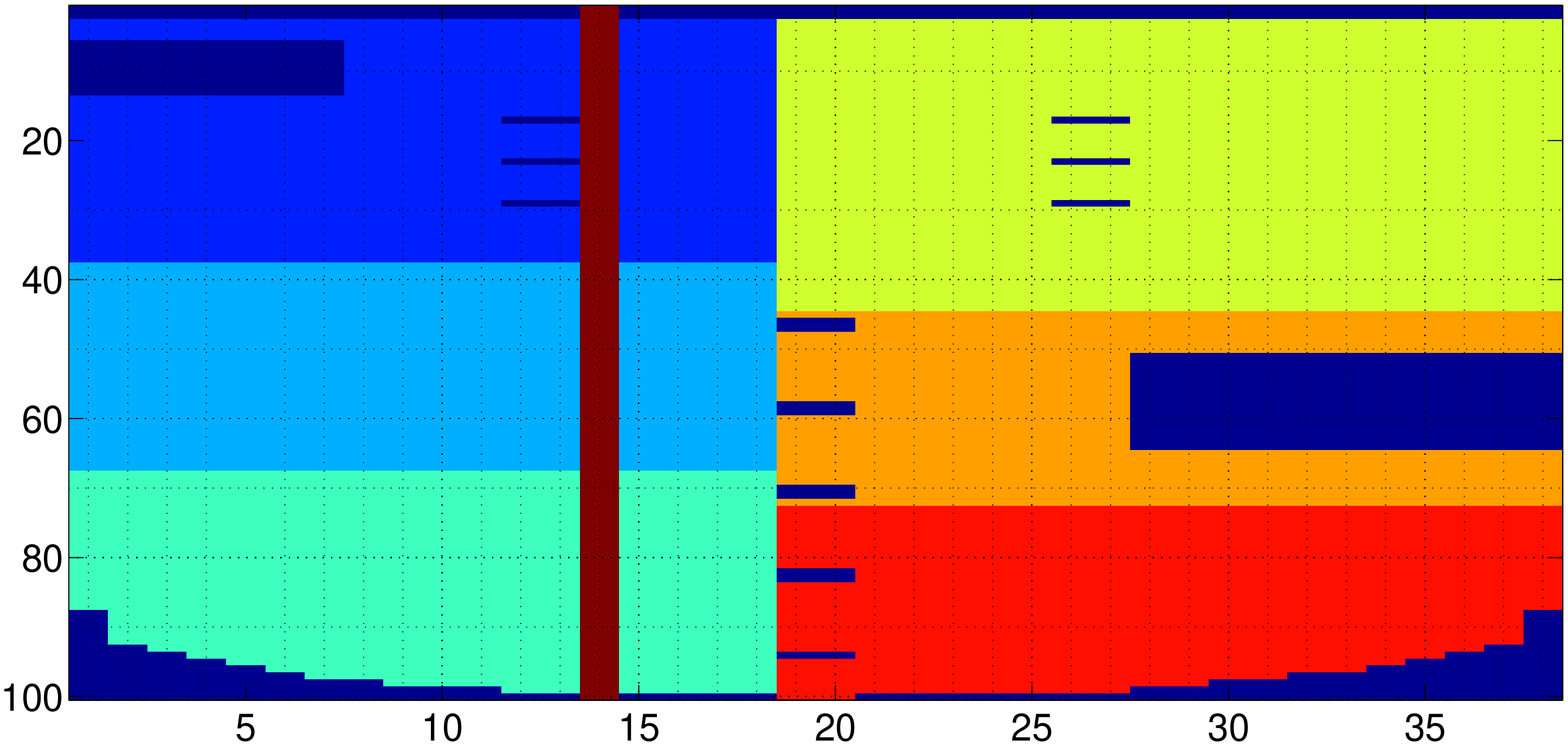}} 
    \subfigure[Problem 4]{\includegraphics[width=0.49\textwidth]{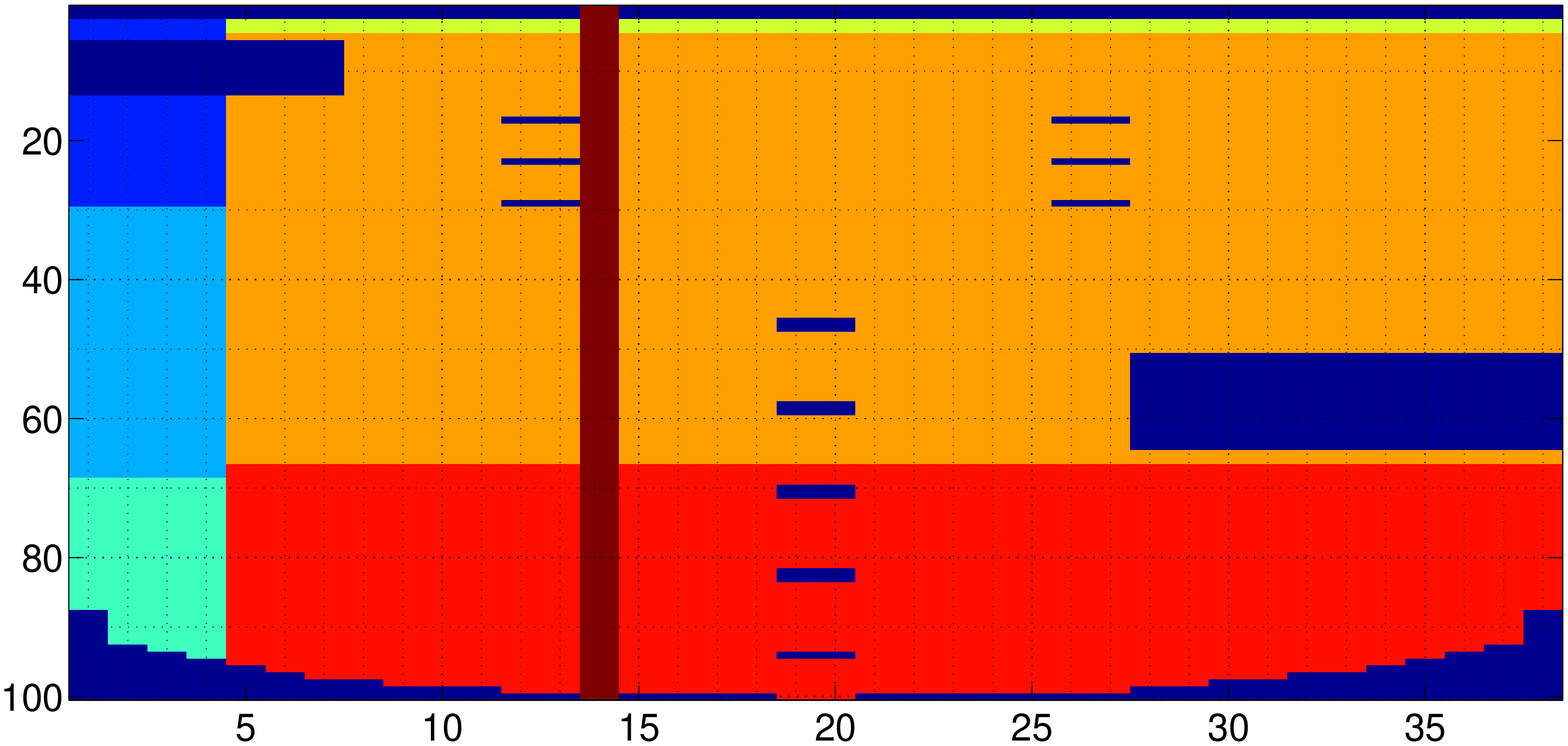}}
    \subfigure[Problem 5]{\includegraphics[width=0.49\textwidth]{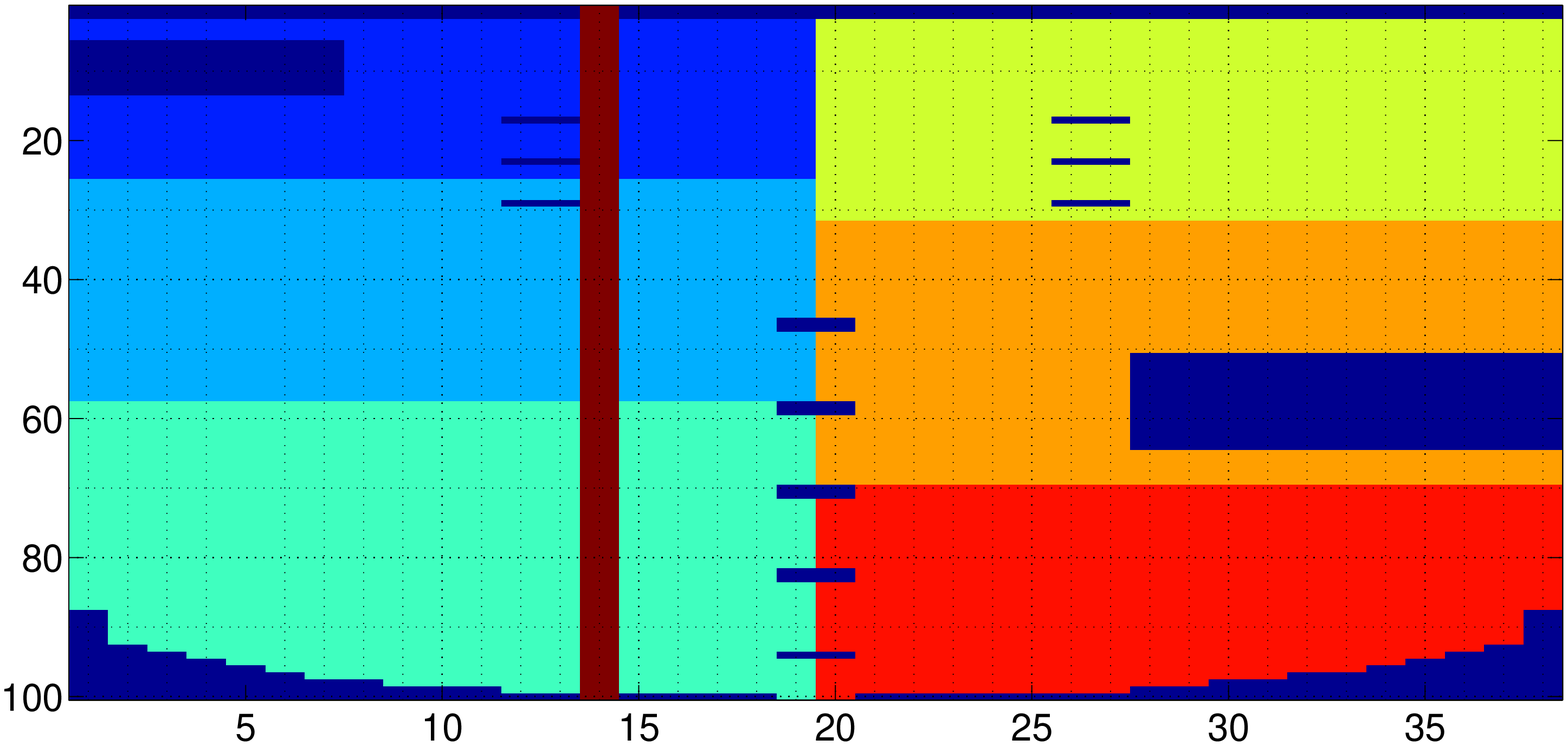}} 
    \subfigure[Problem 6]{\includegraphics[width=0.49\textwidth]{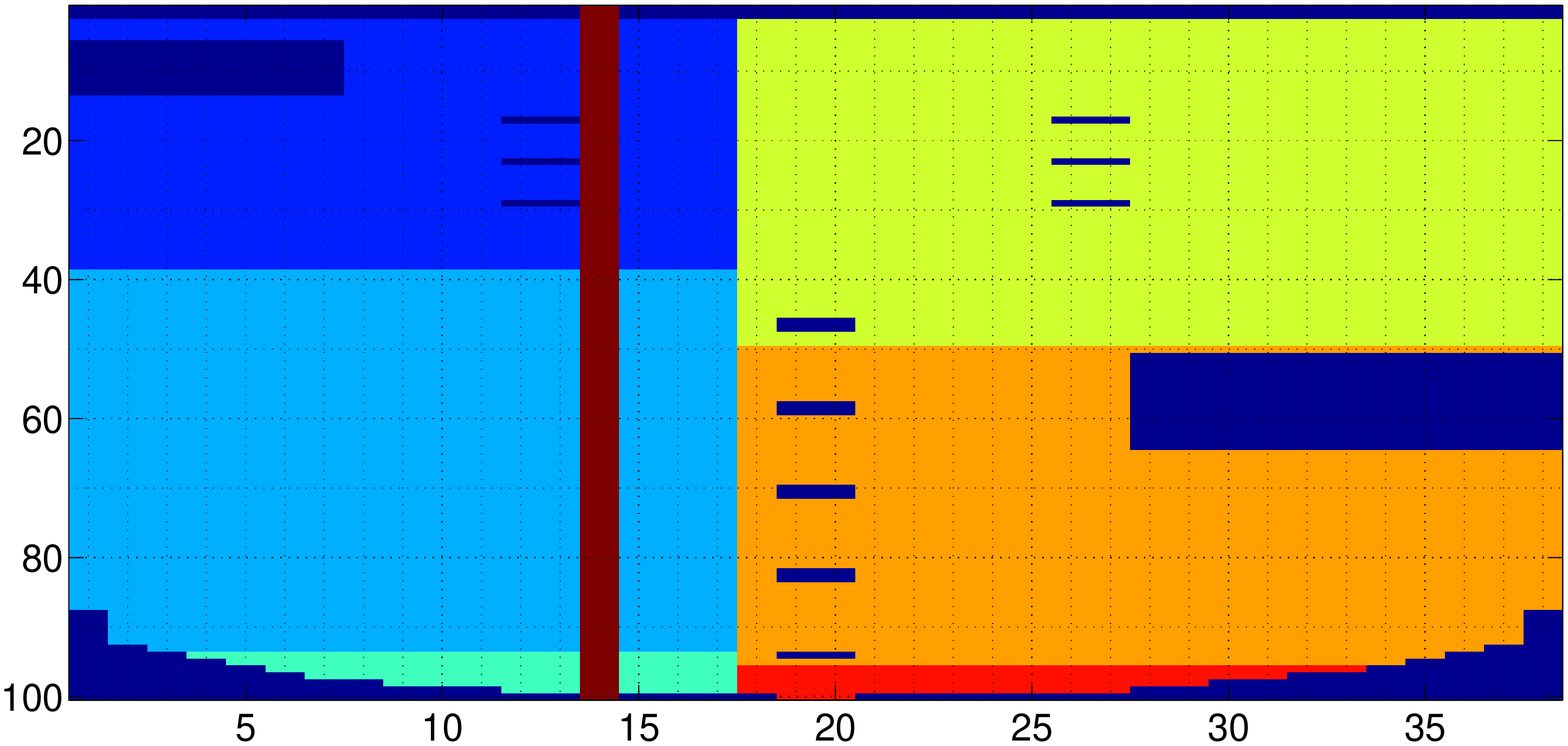}} 
    \caption{Optimal Solutions problems 1-6}
    \label{fig:foobar}
\end{figure}

\begin{figure}
    \centering
    \subfigure[$10^{th}$ Iteration-Prob.1]{\includegraphics[width=0.49\textwidth]{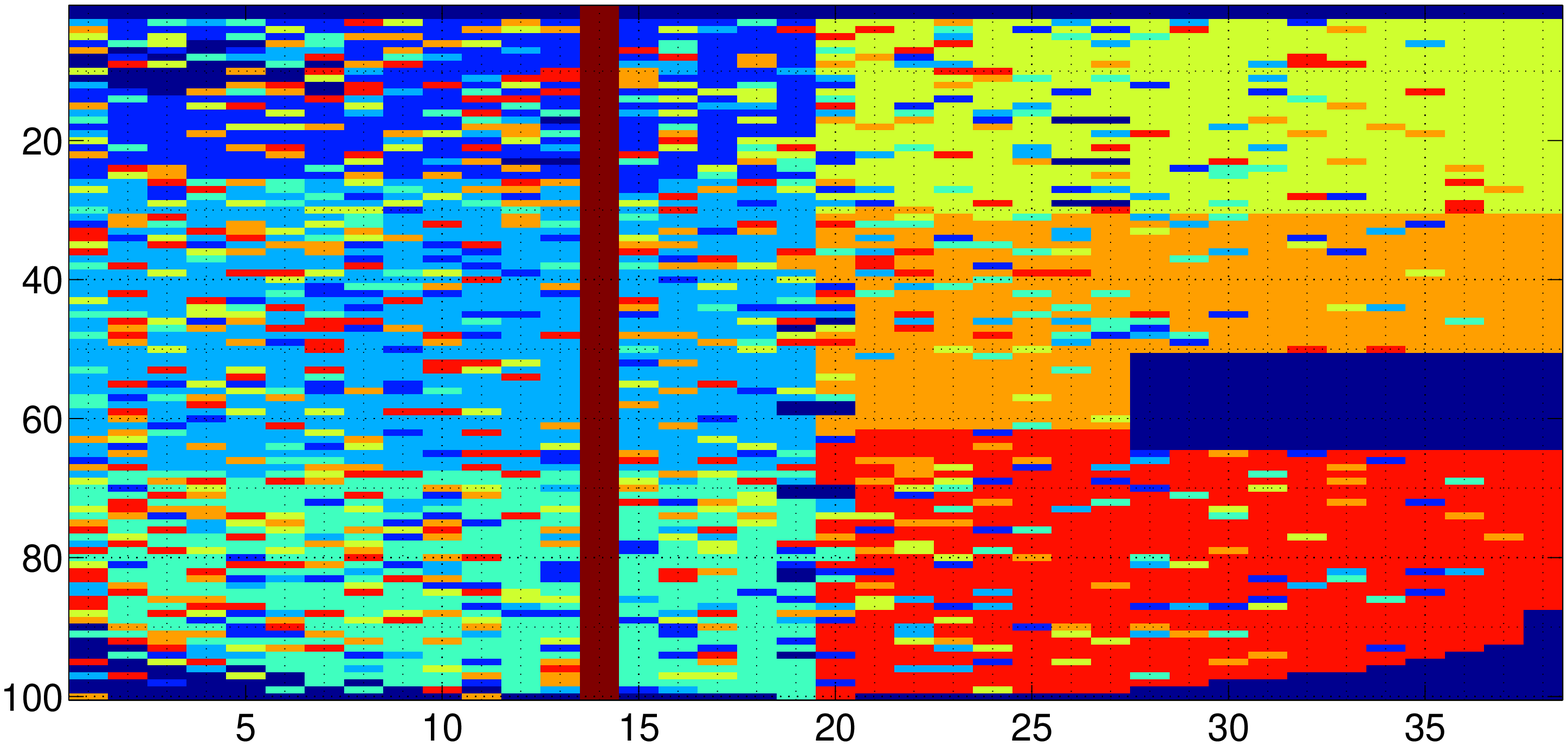}} 
    \subfigure[Best Solution- Prob.1]{\includegraphics[width=0.49\textwidth]{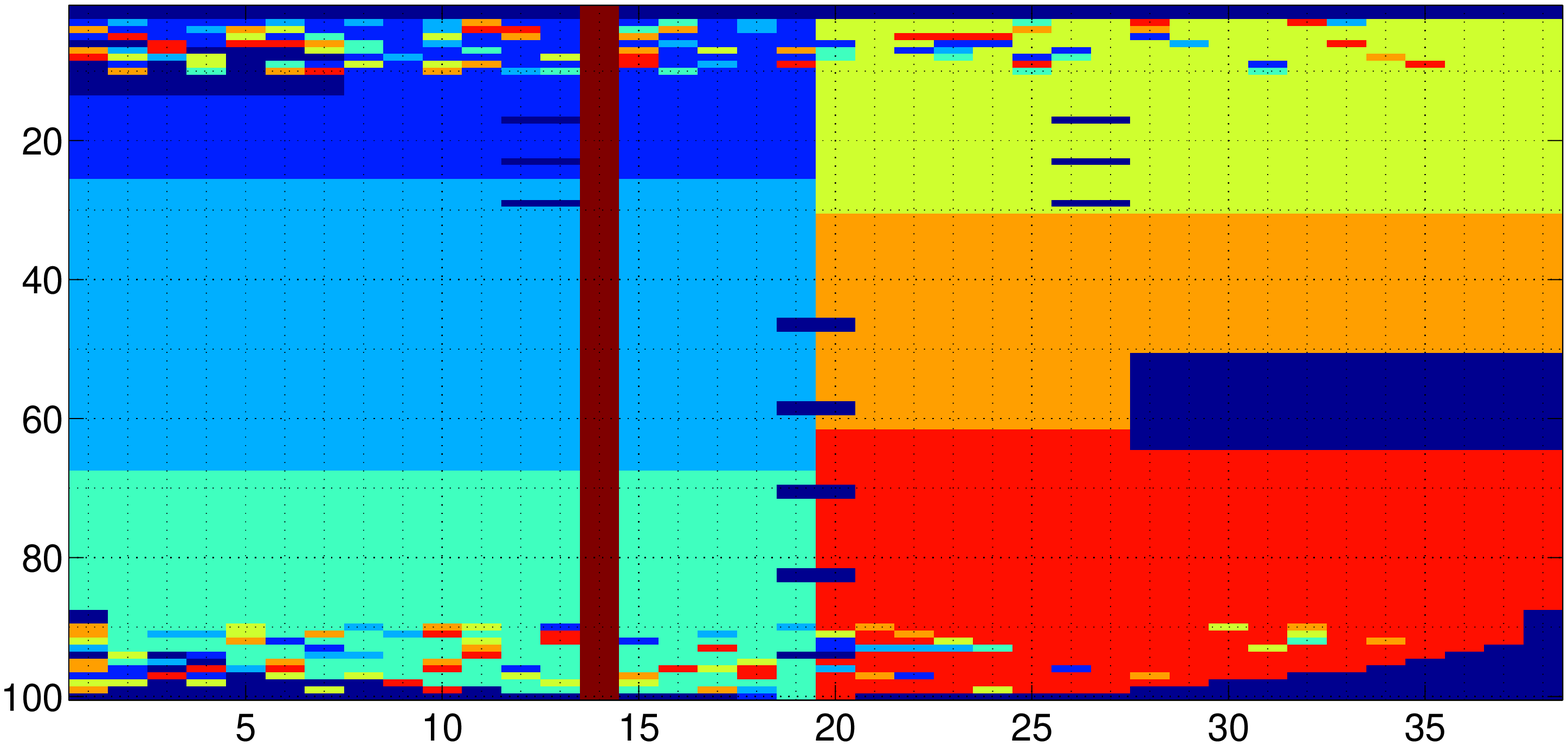}}
    \subfigure[$10^{th}$ Iteration-Prob.2]{\includegraphics[width=0.49\textwidth]{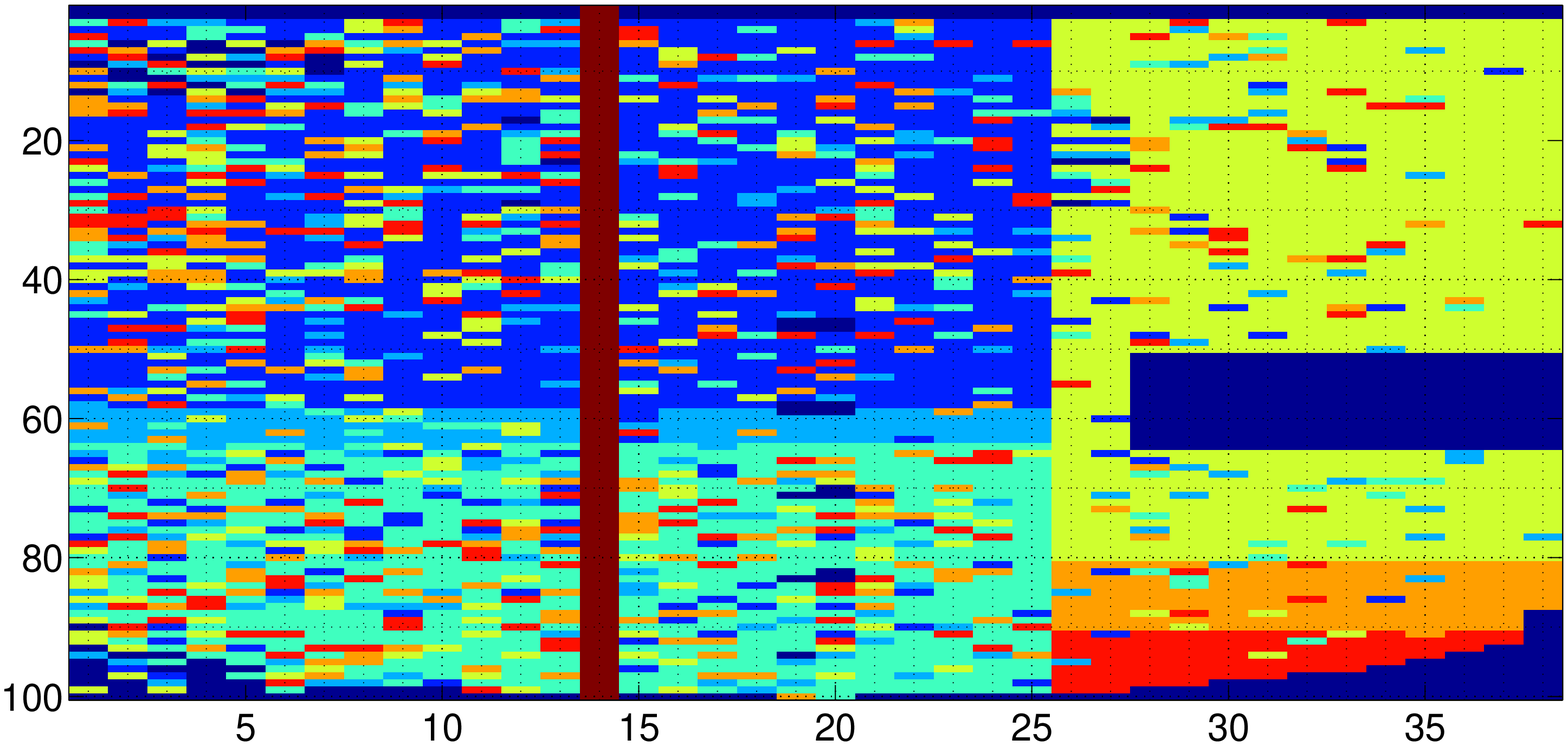}} 
    \subfigure[Best Solution- Prob.2]{\includegraphics[width=0.49\textwidth]{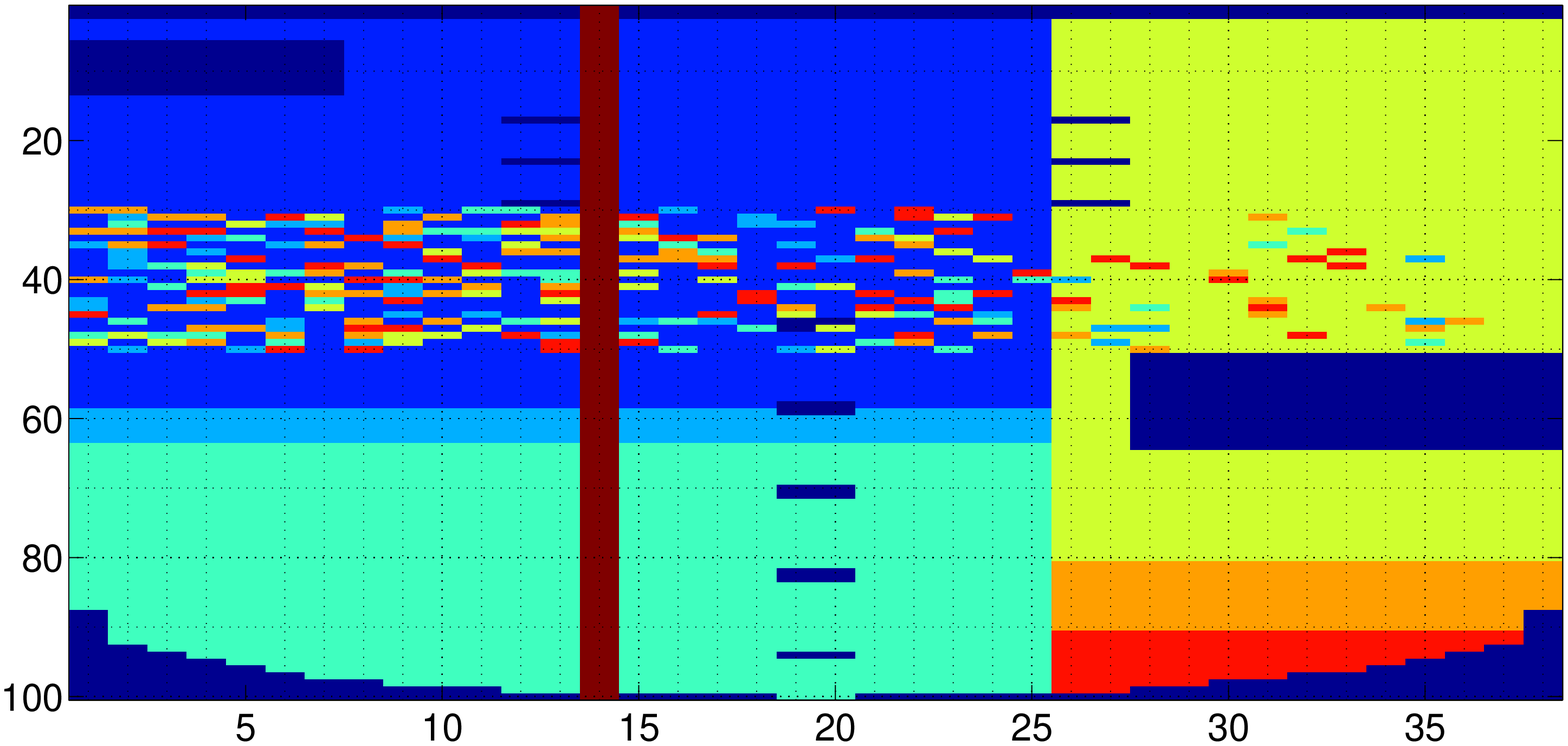}}
    \subfigure[$10^{th}$ Iteration-Prob.3]{\includegraphics[width=0.49\textwidth]{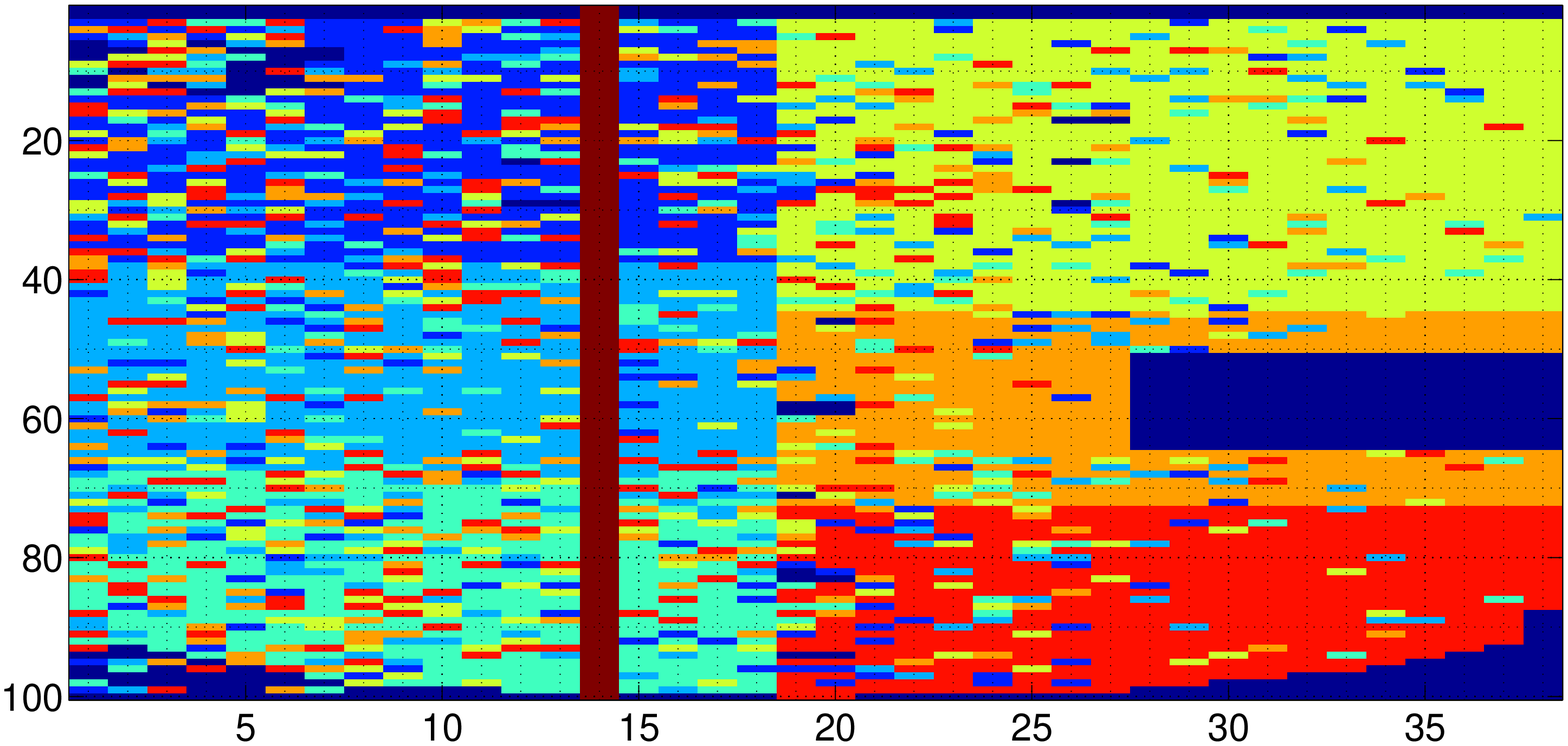}} 
    \subfigure[Best Solution- Prob.3]{\includegraphics[width=0.49\textwidth]{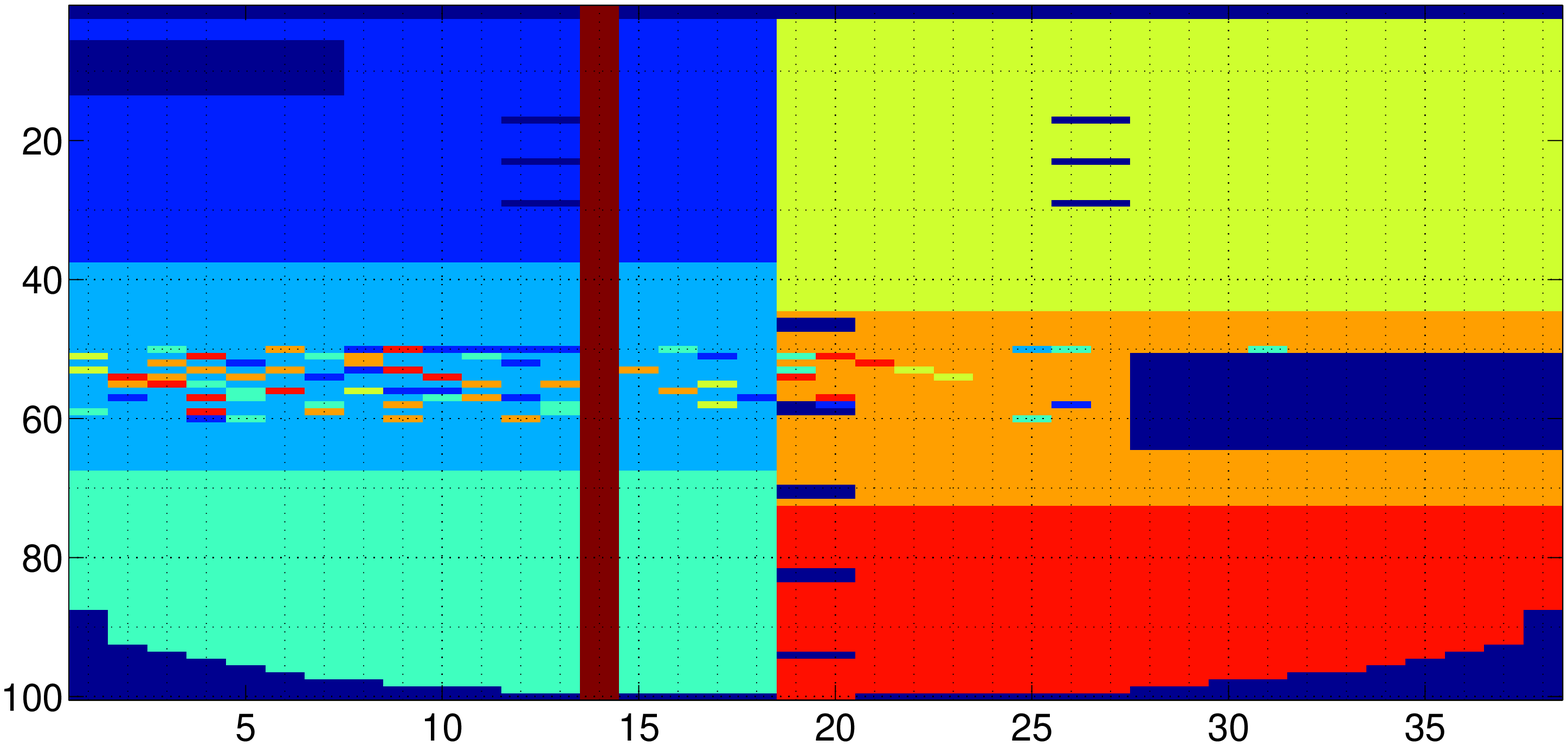}}
    \subfigure[$10^{th}$ Iteration-Prob.4]{\includegraphics[width=0.49\textwidth]{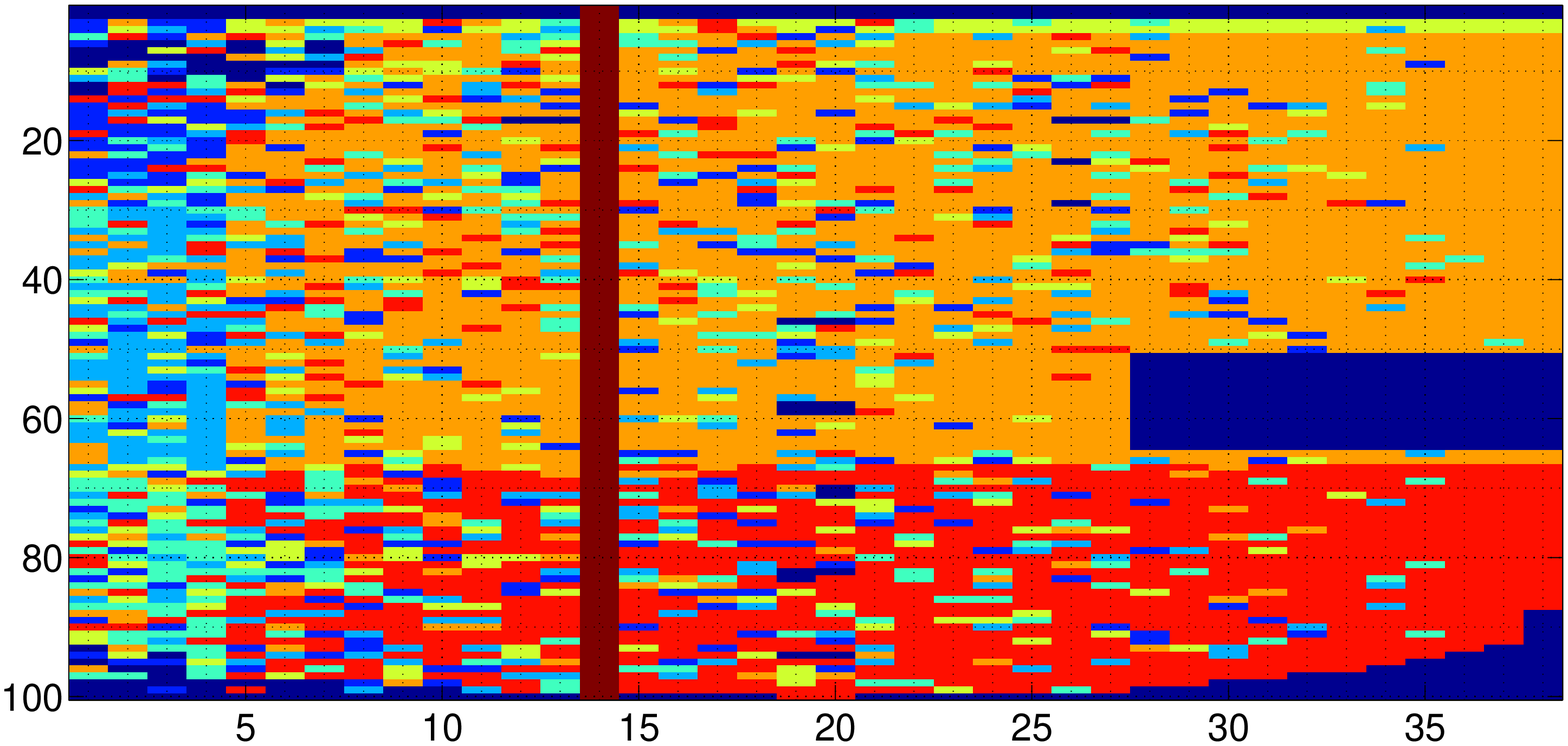}} 
    \subfigure[Best Solution- Prob.4]{\includegraphics[width=0.49\textwidth]{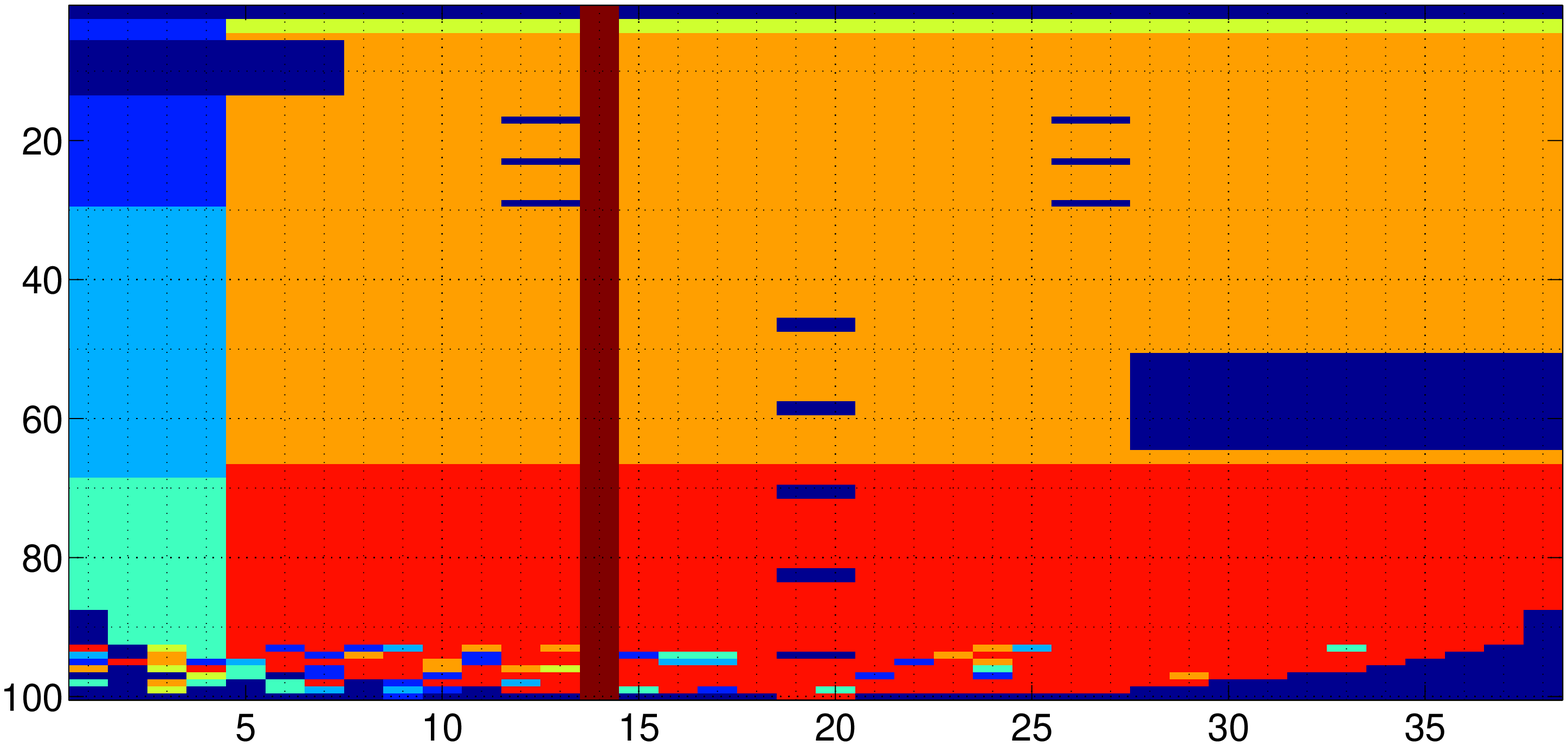}}
    \subfigure[$10^{th}$ Iteration-Prob.5]{\includegraphics[width=0.49\textwidth]{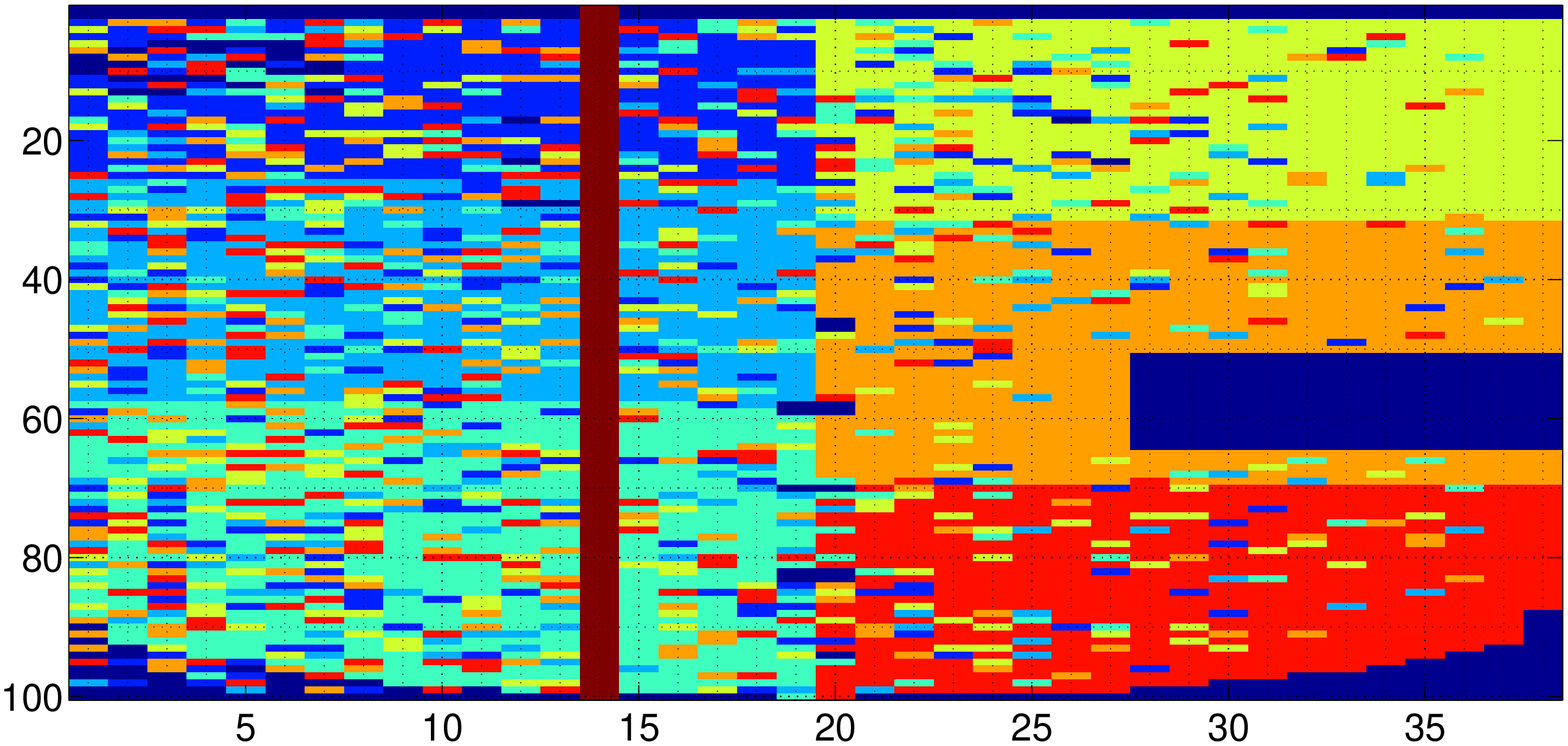}} 
    \subfigure[Best Solution- Prob.5]{\includegraphics[width=0.49\textwidth]{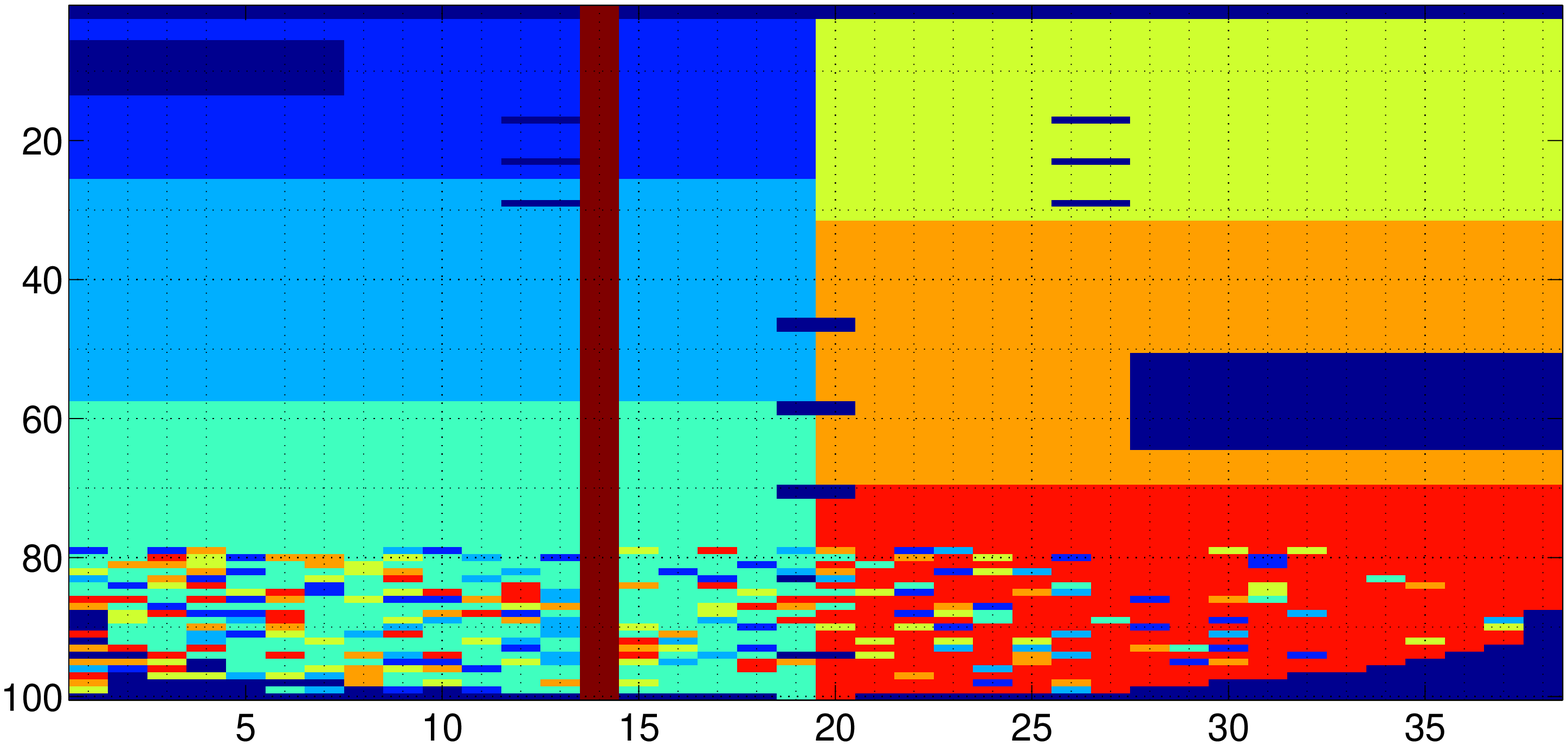}}
    \subfigure[$10^{th}$ Iteration-Prob.6]{\includegraphics[width=0.49\textwidth]{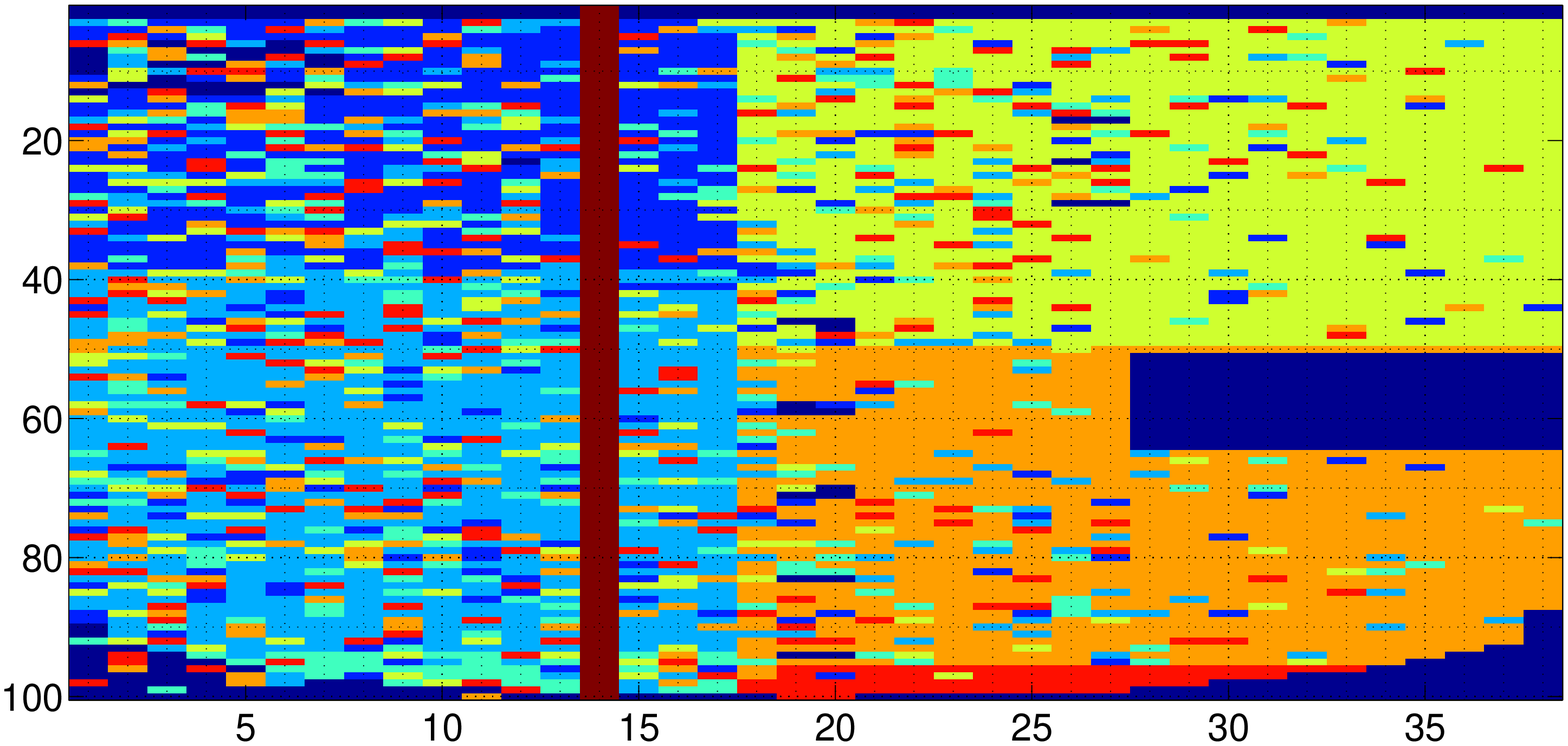}} 
    \subfigure[Best Solution- Prob.6]{\includegraphics[width=0.49\textwidth]{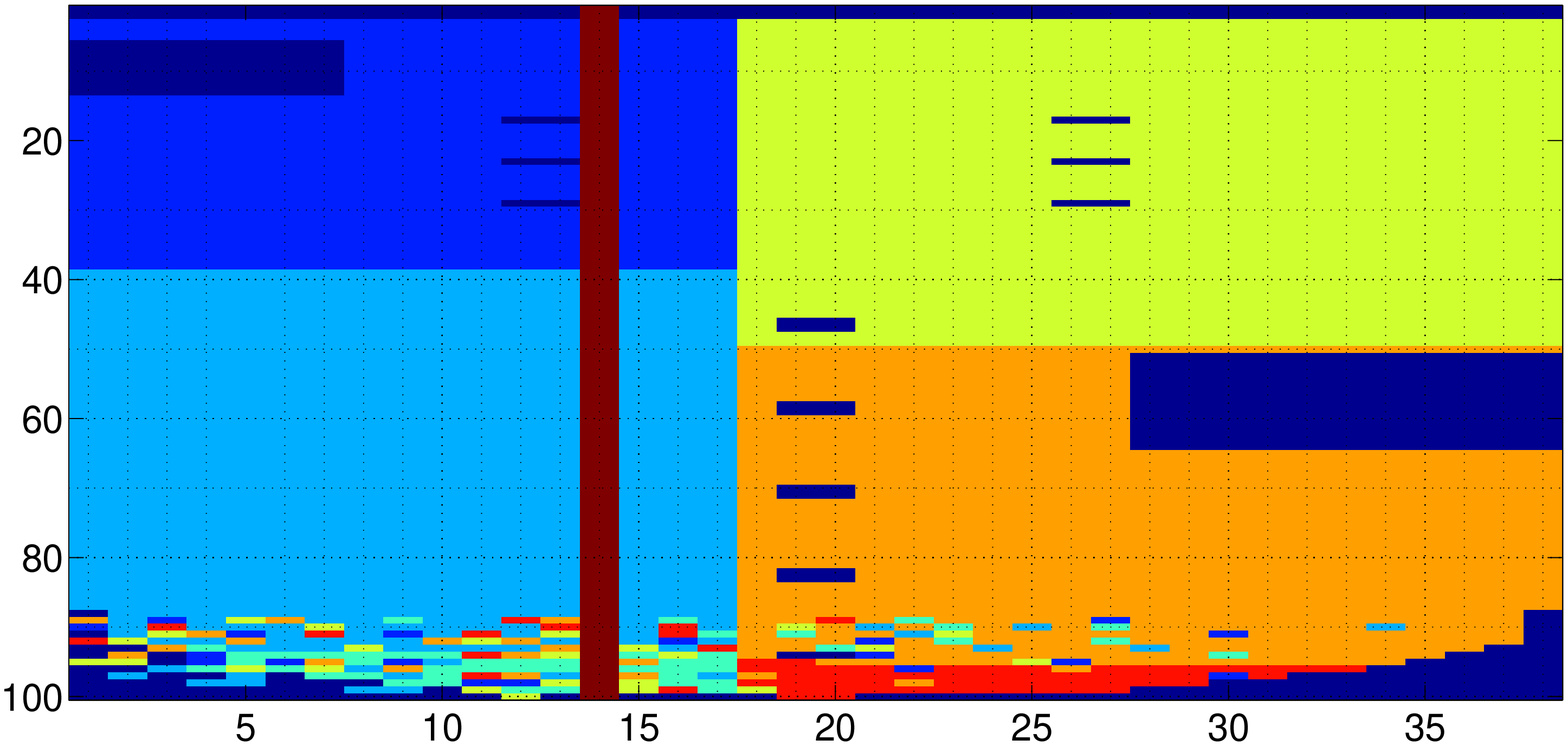}}
    \caption{Behavior of LCA for benchmarks}
    \label{fig:foobar}
\end{figure}

\begin{figure}
    \centering
    \subfigure[$10^{th}$ Iteration-Prob.1]{\includegraphics[width=0.49\textwidth]{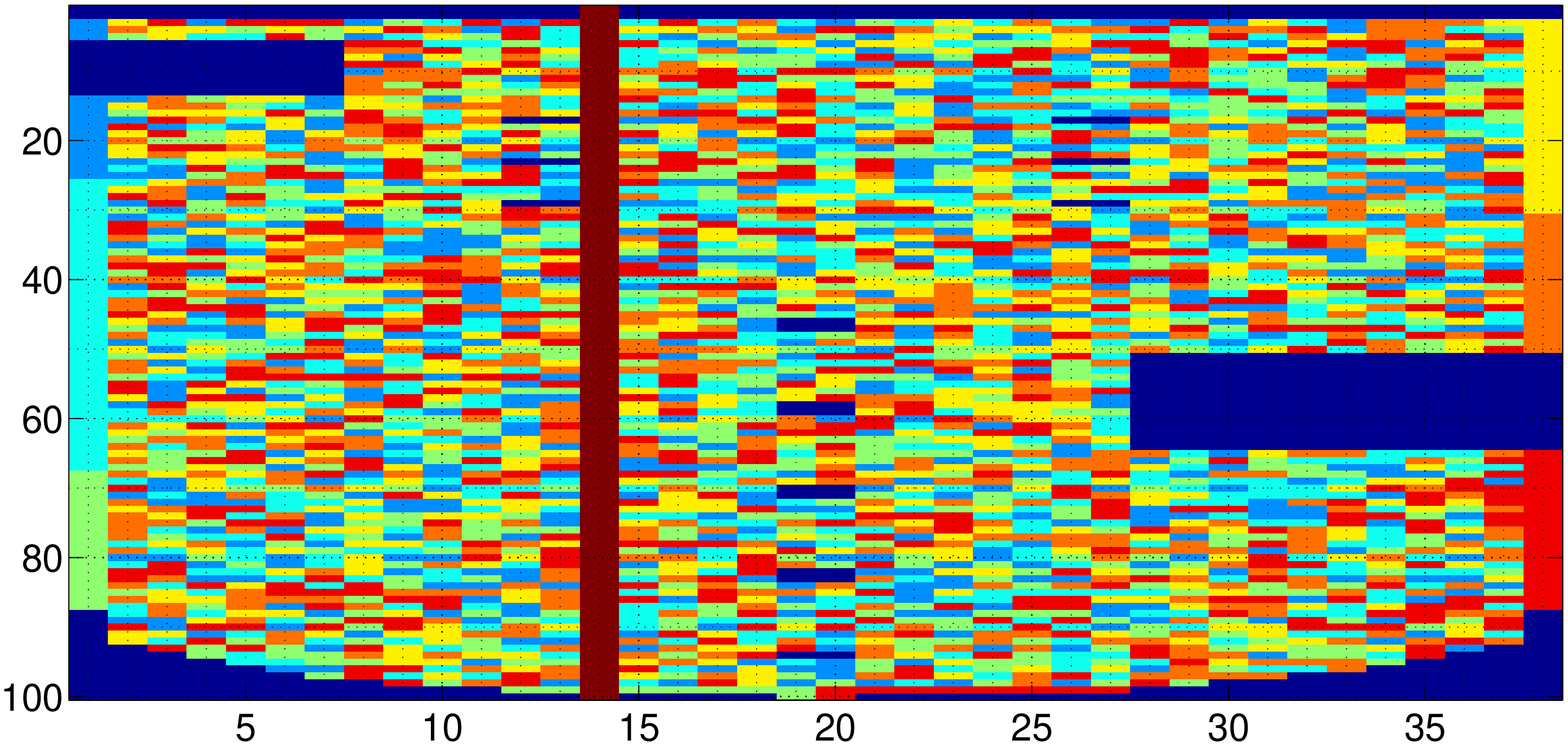}} 
    \subfigure[Best Solution- Prob.1]{\includegraphics[width=0.49\textwidth]{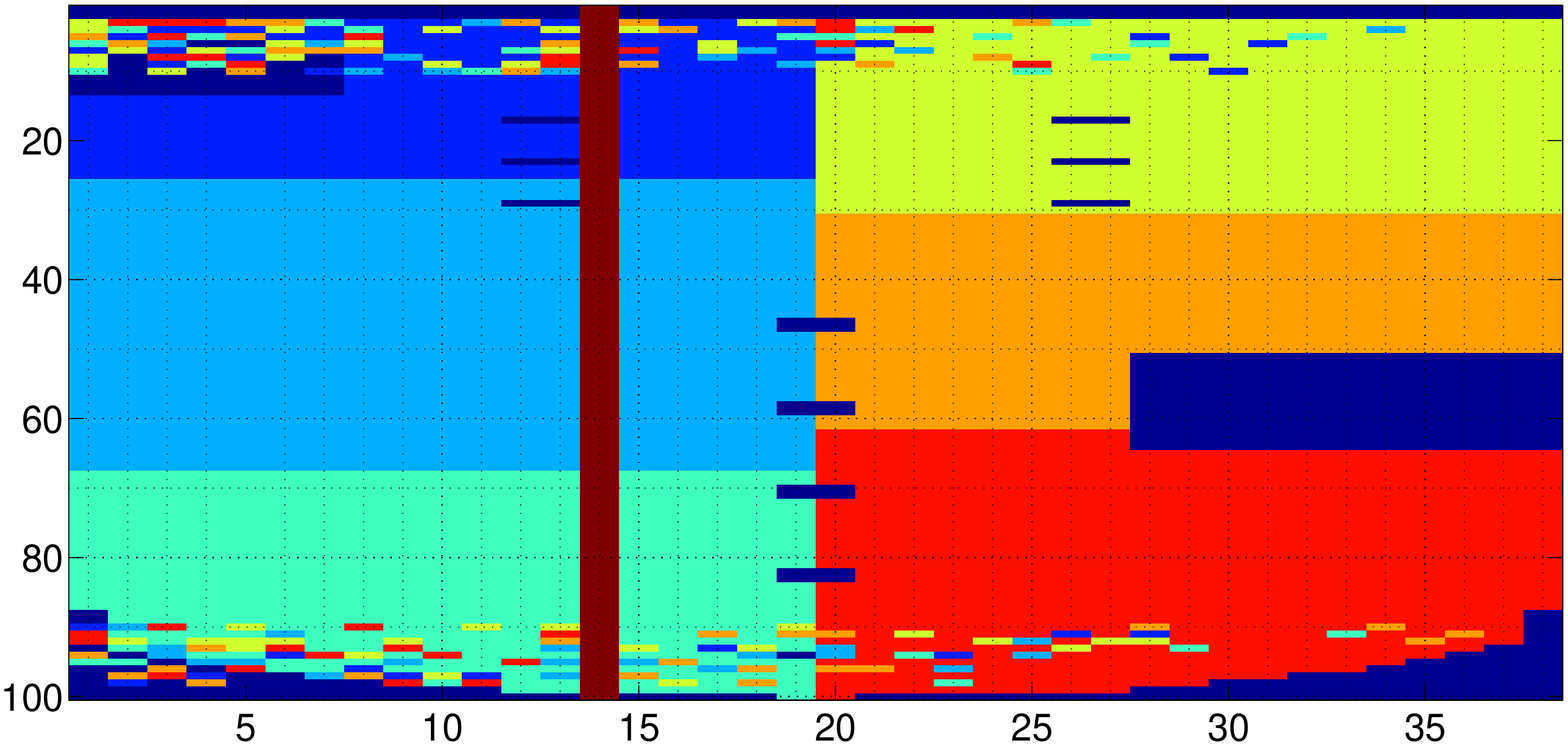}}
    \subfigure[$10^{th}$ Iteration-Prob.2]{\includegraphics[width=0.49\textwidth]{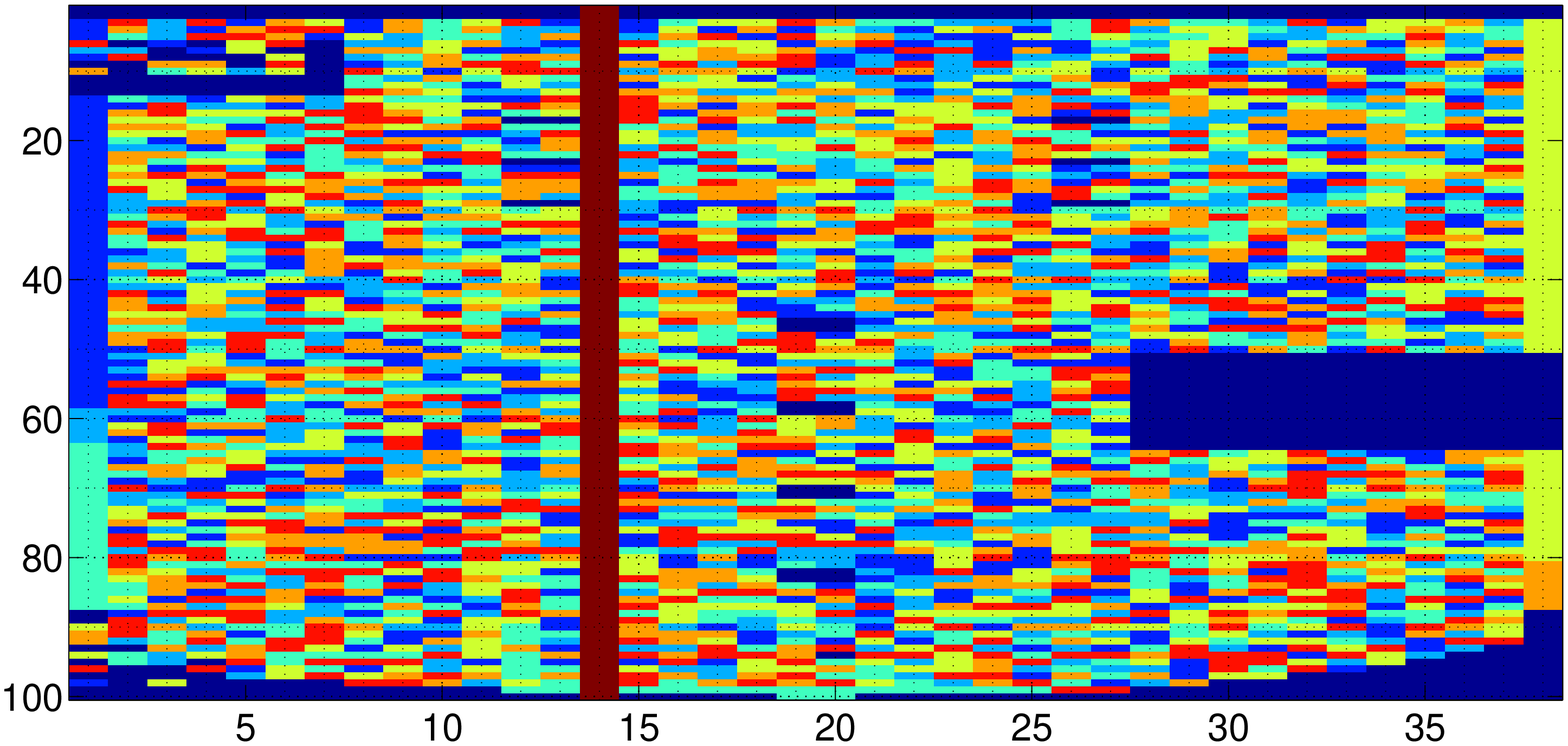}} 
    \subfigure[Best Solution- Prob.2]{\includegraphics[width=0.49\textwidth]{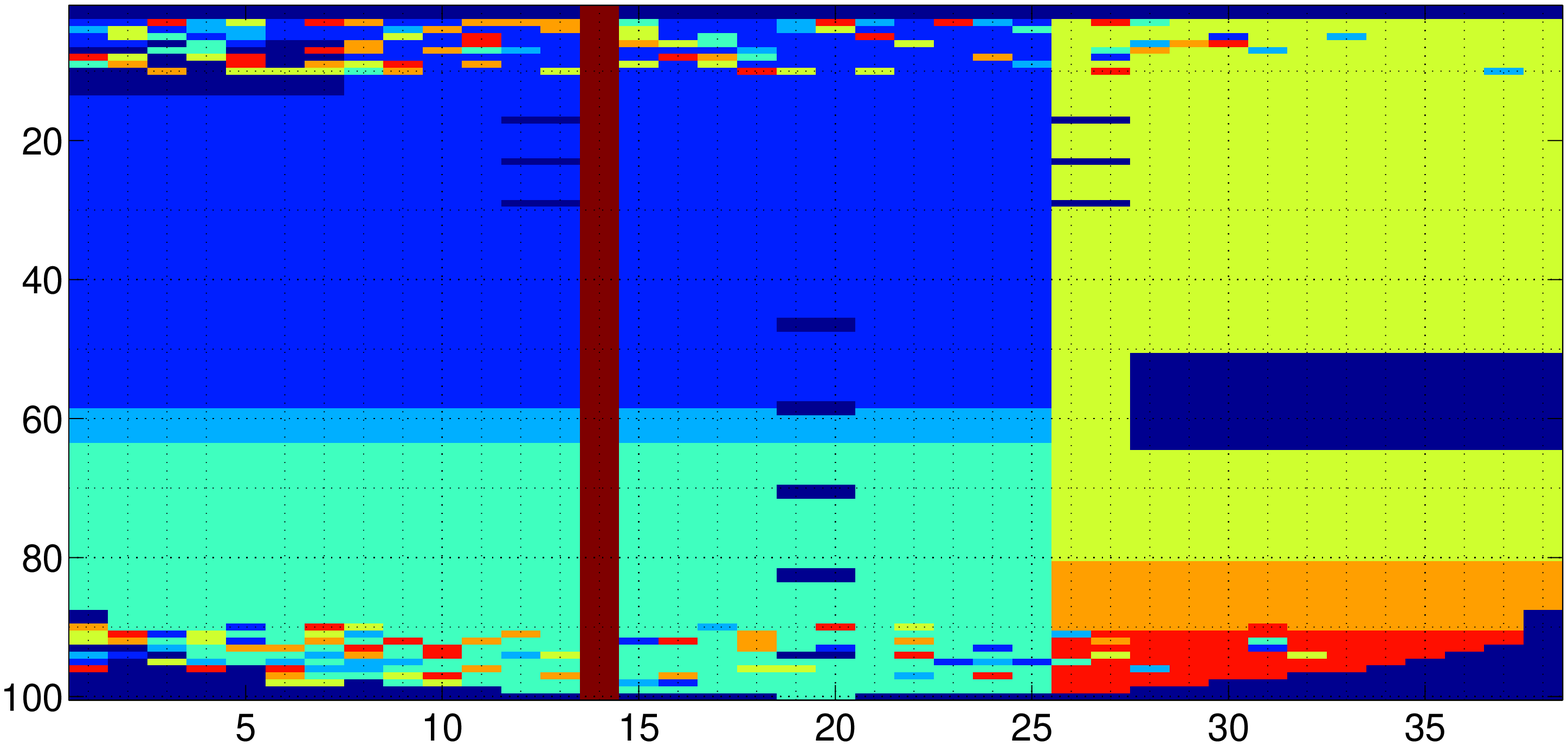}}
    \subfigure[$10^{th}$ Iteration-Prob.3]{\includegraphics[width=0.49\textwidth]{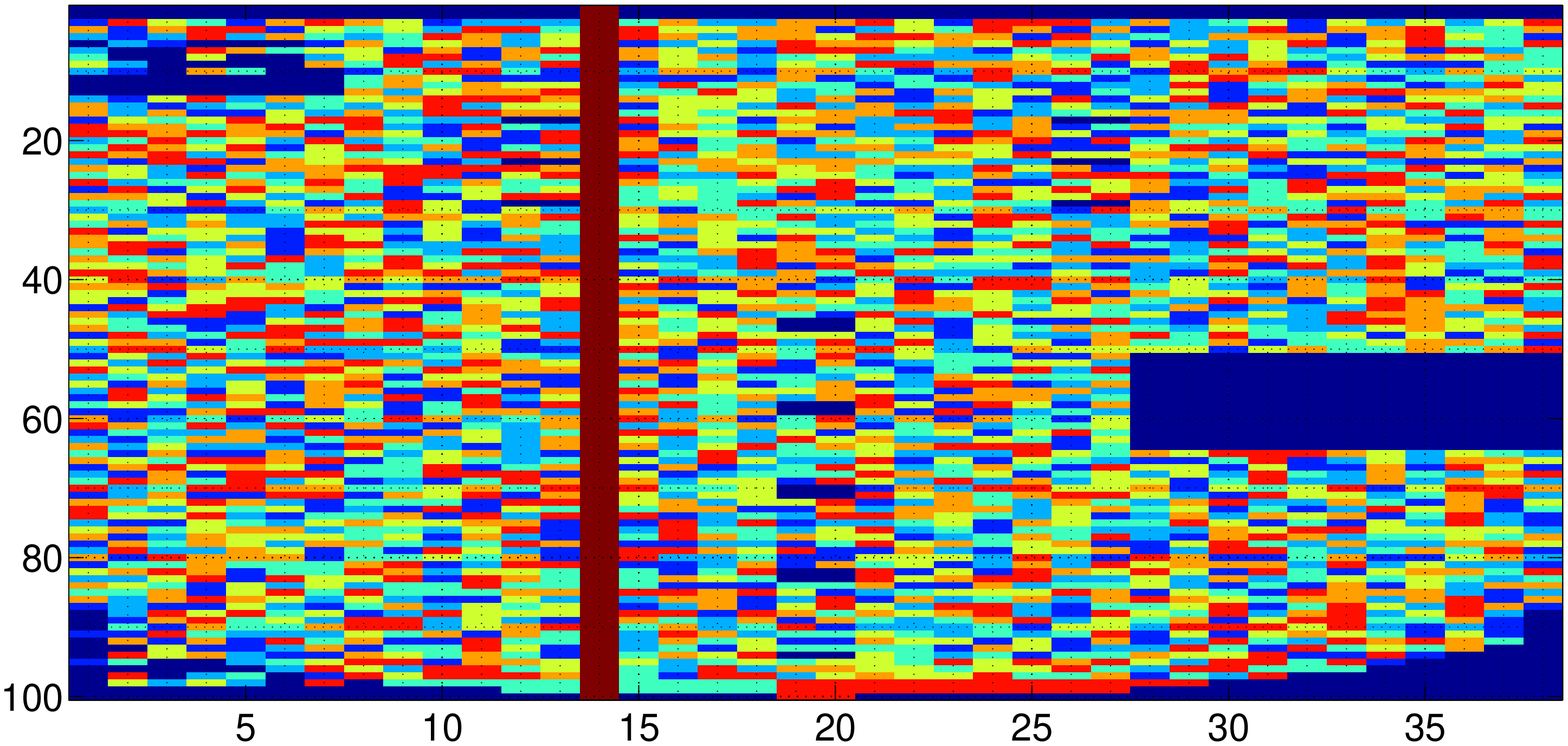}} 
    \subfigure[Best Solution- Prob.3]{\includegraphics[width=0.49\textwidth]{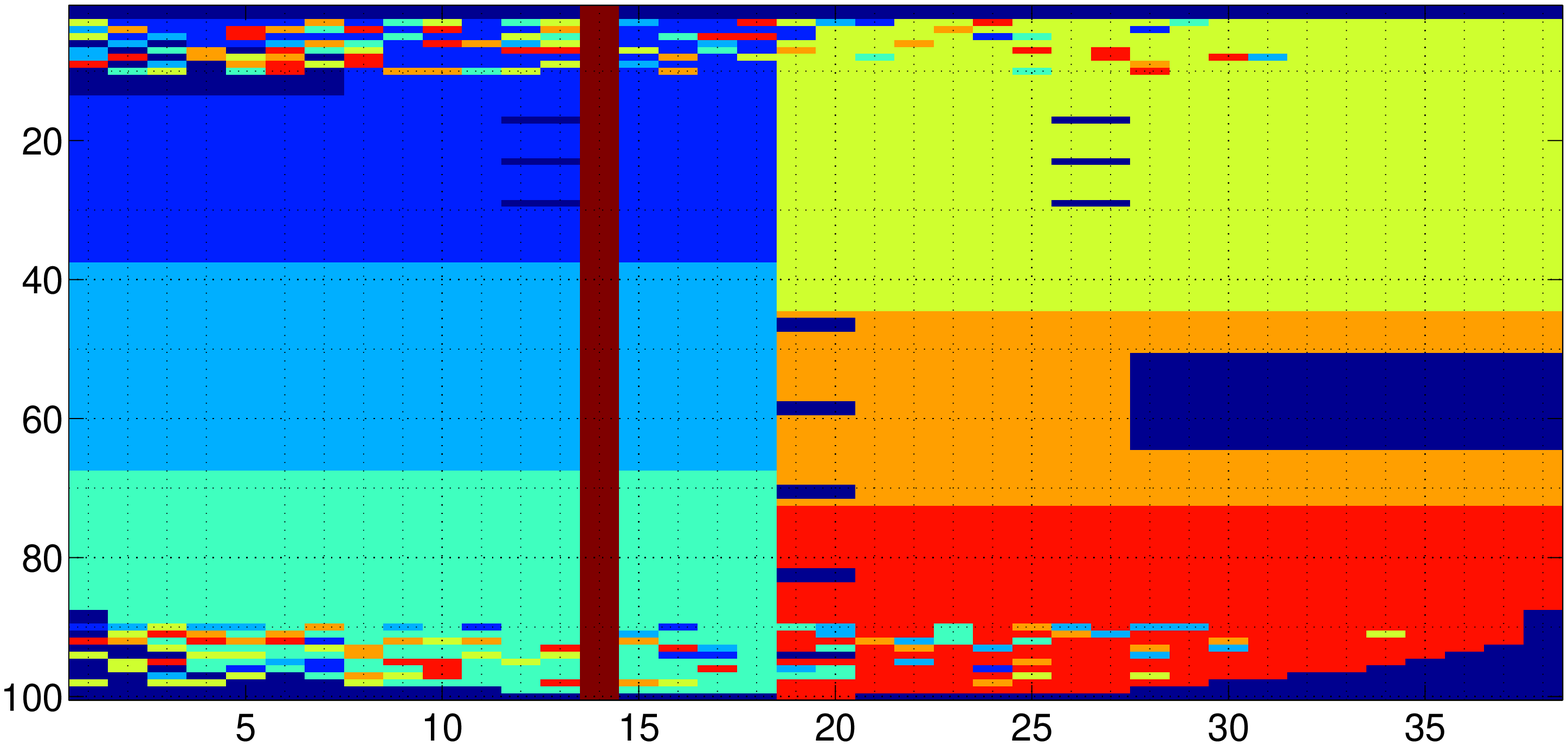}}
    \subfigure[$10^{th}$ Iteration-Prob.4]{\includegraphics[width=0.49\textwidth]{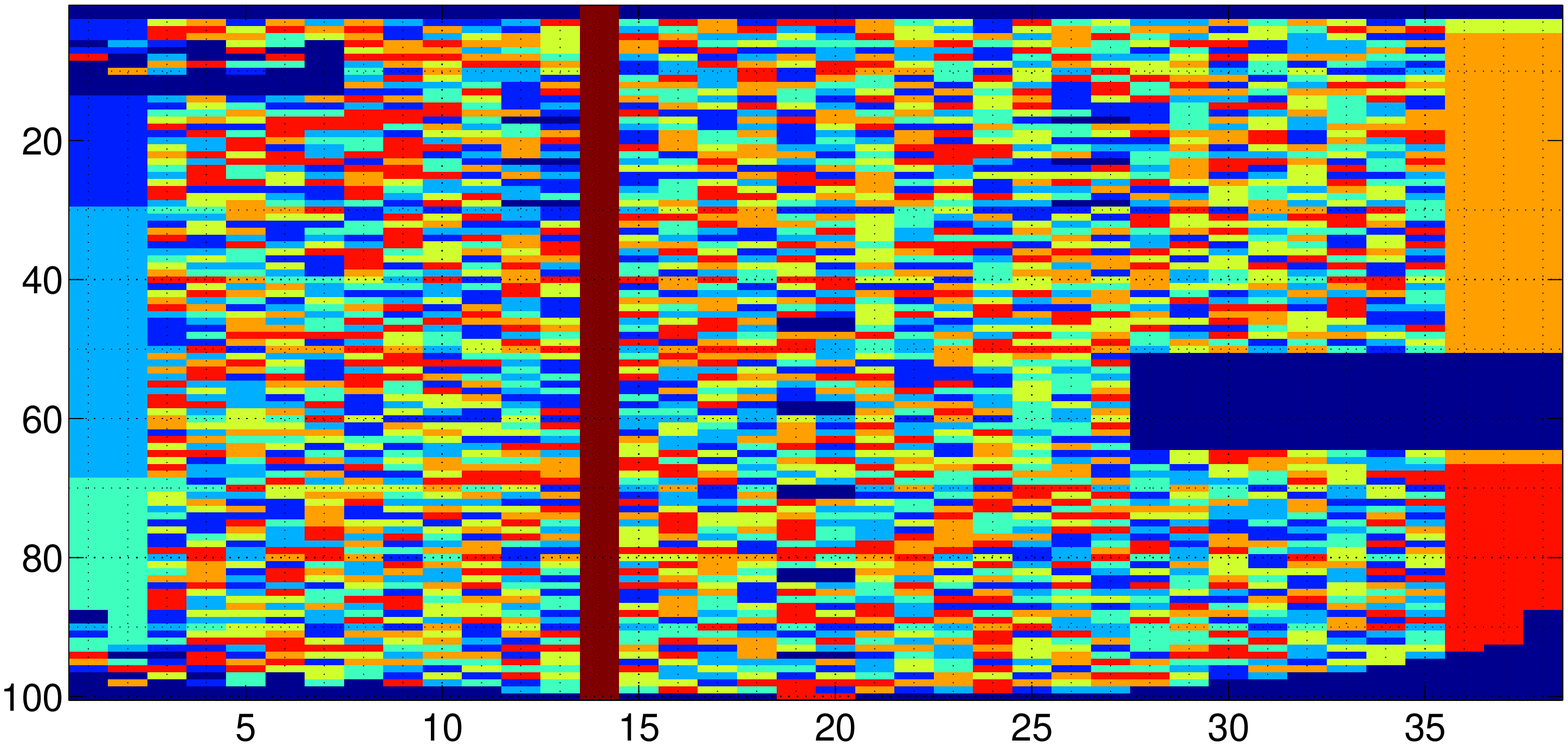}} 
    \subfigure[Best Solution- Prob.4]{\includegraphics[width=0.49\textwidth]{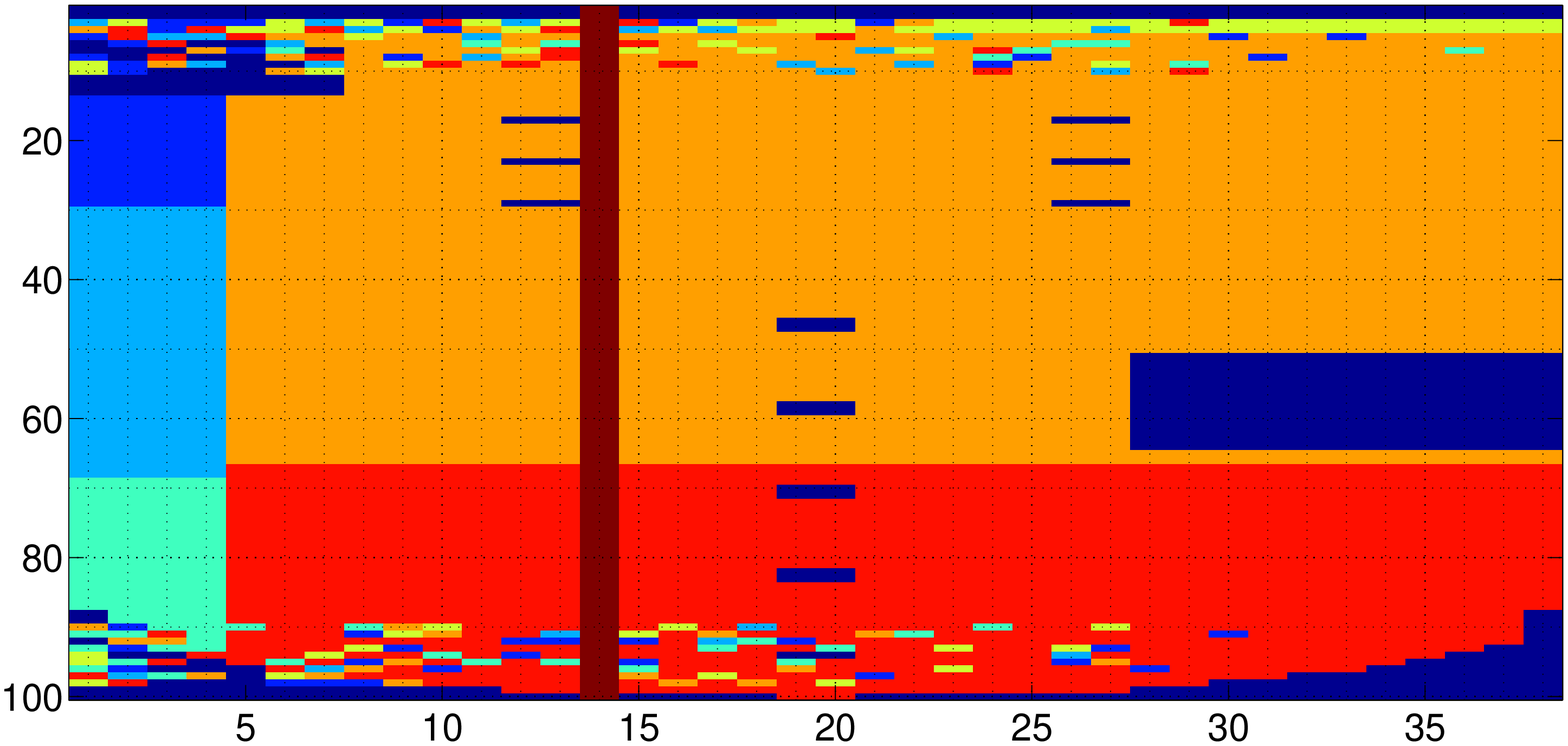}}
    \subfigure[$10^{th}$ Iteration-Prob.5]{\includegraphics[width=0.49\textwidth]{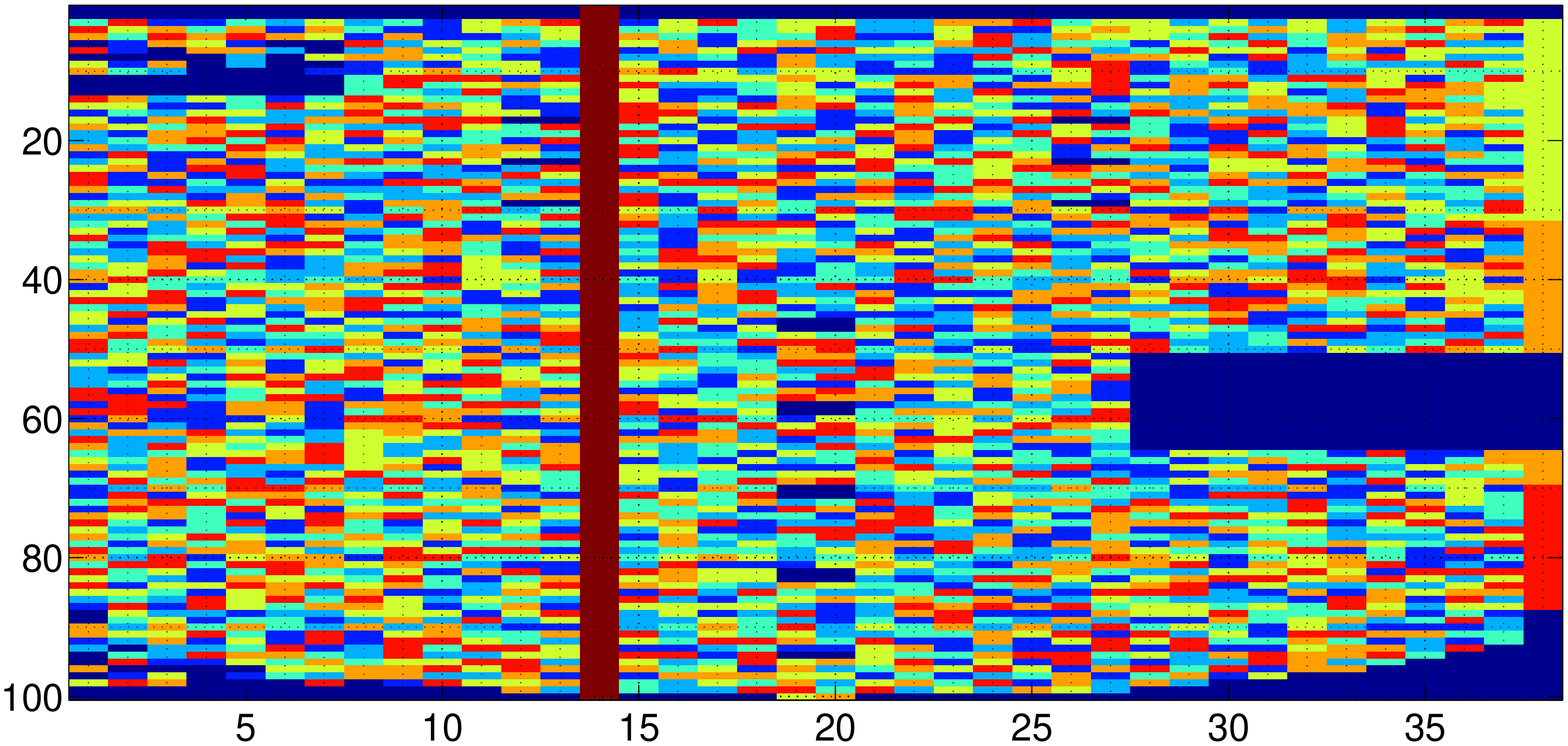}} 
    \subfigure[Best Solution- Prob.5]{\includegraphics[width=0.49\textwidth]{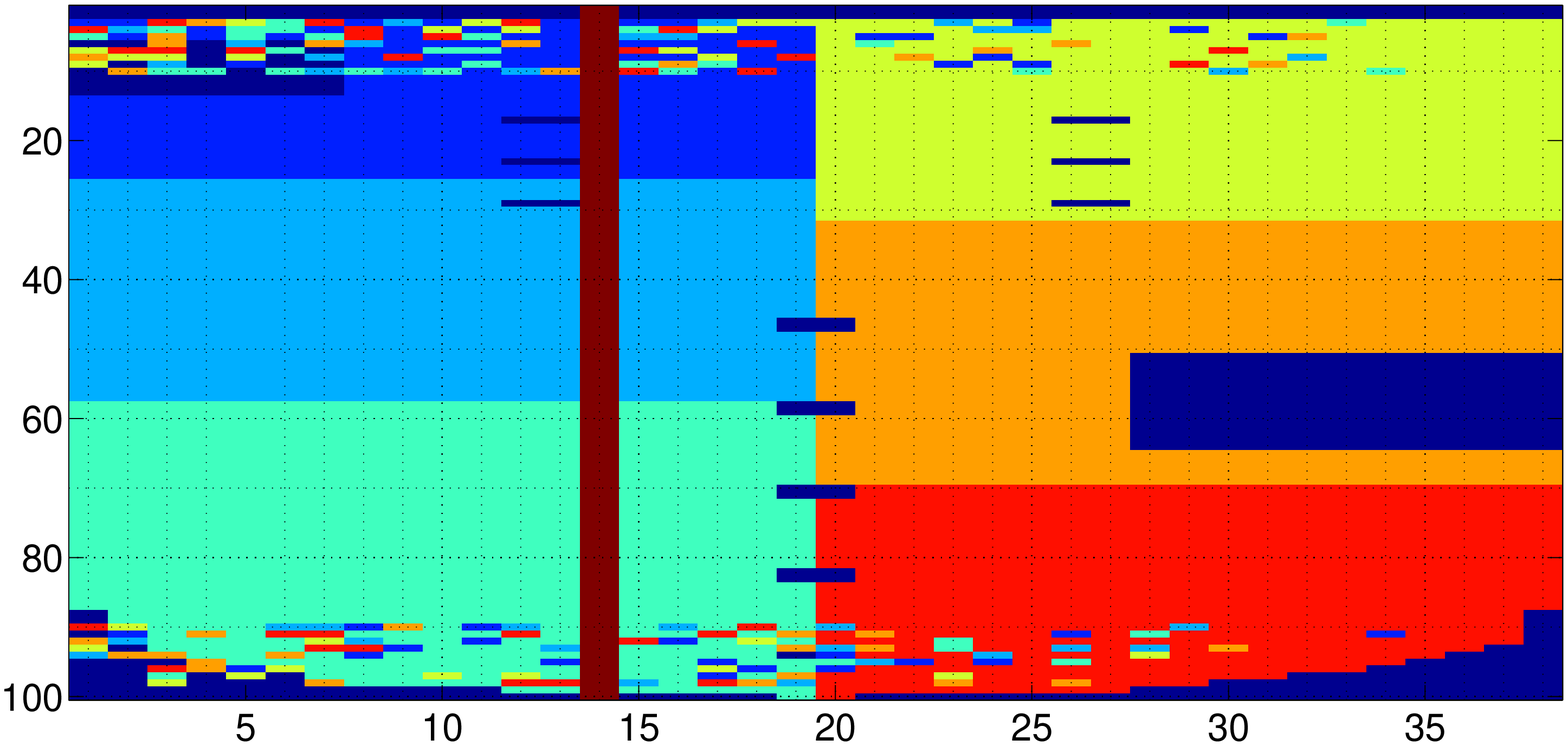}}
    \subfigure[$10^{th}$ Iteration-Prob.6]{\includegraphics[width=0.49\textwidth]{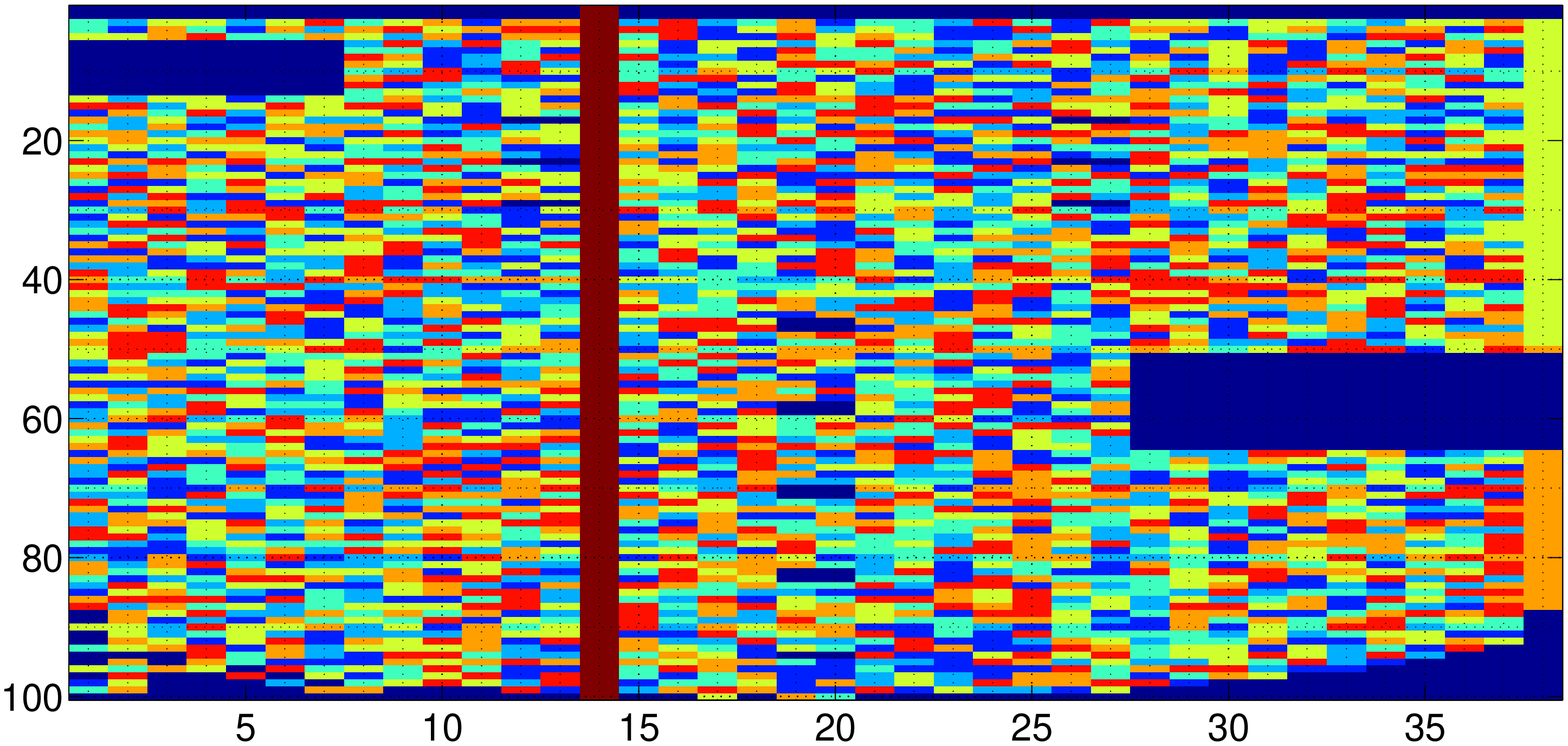}} 
    \subfigure[Best Solution- Prob.6]{\includegraphics[width=0.49\textwidth]{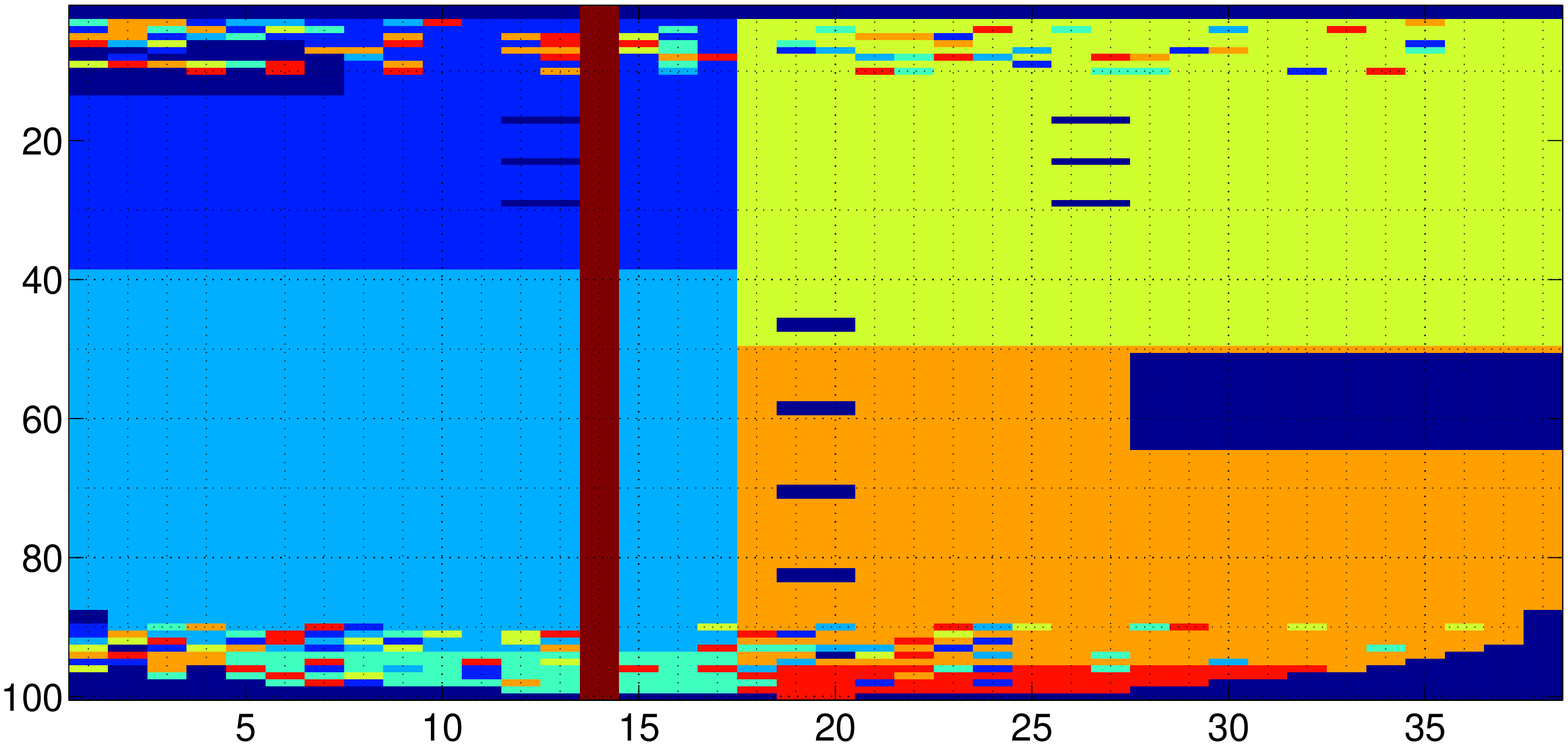}}
    \caption{Behavior of VEA for benchmarks}
    \label{fig:foobar}
\end{figure}

 \begin{table*} [t]
  	\caption{Results of SA for Benchmarks }
  	\begin{center}
  		\label{relatedwork}
  		\small{ 
  			\begin{tabular}{| p{.18\textwidth} | p{.22\textwidth} | p{.18\textwidth} | p{.18\textwidth} |  p{.18\textwidth} | p{.18\textwidth} | p{.09\textwidth} | p{.09\textwidth} |} 
  				\hline
  				
  				\textbf{Problems} &\textbf{ Initial Solution}  &\textbf{ 5 Iterations }&  \textbf{10 Iterations } & \textbf{Best Solution}\\ 
\hline
  						
  				  Inst0&31650&31350&10360&2570\\ 
\hline  
  				
  				 Inst1&27440&27730&11230&1980\\ 
\hline      
  				
  				Inst2&27970&27980&9610&2380\\  
\hline    

  				 Inst3&27780&27220&12760&3450\\  
\hline    
  				 Inst4&27930&26540&9780&1210\\  
\hline    
  				 Inst5&27790&27990&9110&2340\\  
\hline 

  			\end{tabular}
  			
  		}
  		
  	\end{center}                                                
  	
  \end{table*}

 \begin{table*} [t]
  	\caption{Results of TS for Benchmarks }
  	\begin{center}
  		\label{relatedwork}
  		\small{ 
  			\begin{tabular}{| p{.18\textwidth} | p{.22\textwidth} | p{.18\textwidth} | p{.18\textwidth} |  p{.18\textwidth} | p{.18\textwidth} | p{.09\textwidth} | p{.09\textwidth} |} 
  				\hline
  				
  				\textbf{Problems} &\textbf{ Initial Solution}  &\textbf{ 5 Iterations }&  \textbf{10 Iterations } & \textbf{Best Solution}\\ 
\hline
  						
  				  Inst0&31570&31570&31260&1030\\ 
\hline  
  				
  				 Inst1&31630&31630&31260&1210\\ 
\hline      
  				
  				Inst2&31640&31600&31220&1920\\  
\hline    

  				 Inst3&32070&31880&31210&1880\\  
\hline    
  				 Inst4&32020&31920&31230&2100\\  
\hline    
  				 Inst5&31770&31420&31120&1750\\  
\hline 

  			\end{tabular}
  			
  		}
  		
  	\end{center}                                                
  	
  \end{table*}

 \begin{table*} [t]
  	\caption{Results of LCA for Benchmarks }
  	\begin{center}
  		\label{relatedwork}
  		\small{ 
  			\begin{tabular}{| p{.18\textwidth} | p{.22\textwidth} | p{.18\textwidth} | p{.18\textwidth} |  p{.18\textwidth} | p{.18\textwidth} | p{.09\textwidth} | p{.09\textwidth} |} 
  				\hline
  				
  				\textbf{Problems} &\textbf{ Initial Solution}  &\textbf{ 5 Iterations }&  \textbf{10 Iterations } & \textbf{Best Solution}\\ 
\hline
  						
  				  Inst0&31710&31310&9690&1740\\ 
\hline  
  				
  				 Inst1&31600&31230&9720&2370\\ 
\hline      
  				
  				Inst2&31720&31290&9760&660\\  
\hline    

  				 Inst3&31180&30950&9750&510\\  
\hline    
  				 Inst4&26830&24120&9630&2060\\  
\hline    
  				 Inst5&27130&26720&9680&930\\  
\hline 
  			\end{tabular}
  			
  		}
  		
  	\end{center}                                                
  	
  \end{table*}

 \begin{table*} [t]
  	\caption{Results of VEA for Benchmarks }
  	\begin{center}
  		\label{relatedwork}
  		\small{ 
  			\begin{tabular}{| p{.18\textwidth} | p{.22\textwidth} | p{.18\textwidth} | p{.18\textwidth} |  p{.18\textwidth} | p{.18\textwidth} | p{.09\textwidth} | p{.09\textwidth} |} 
  				\hline
  				
  				\textbf{Problems} &\textbf{ Initial Solution}  &\textbf{ 5 Iterations }&  \textbf{10 Iterations } & \textbf{Best Solution}\\ 
\hline
  						
  				  Inst0&30590&28710&9690&1620\\ 
\hline  
  				
  				 Inst1&26380&23560&9200&1630\\ 
\hline      
  				
  				Inst2&27420&24670&7420&620\\  
\hline    

  				 Inst3&24330&21450&6740&440\\  
\hline    
  				 Inst4&26690&22680&9390&1730\\  
\hline    
  				 Inst5&26430&25270&9310&680\\  
\hline 

  			\end{tabular}
  			
  		}
  		
  	\end{center}                                                
  	
  \end{table*}

 \begin{table*} [t]
  	\caption{Results of MVA for Benchmarks }
  	\begin{center}
  		\label{relatedwork}
  		\small{ 
  			\begin{tabular}{| p{.18\textwidth} | p{.22\textwidth} | p{.18\textwidth} | p{.18\textwidth} |  p{.18\textwidth} | p{.18\textwidth} | p{.09\textwidth} | p{.09\textwidth} |} 
  				\hline
  				
  				\textbf{Problems} &\textbf{ Initial Solution}  &\textbf{ 5 Iterations }&  \textbf{10 Iterations } & \textbf{Best Solution}\\ 
\hline
  						
  				  Inst0&31900&31510&7420&1300\\ 
\hline  
  				
  				 Inst1&30300&29880&8760&1150\\ 
\hline      
  				
  				Inst2&28720&28210&8560&380\\  
\hline    

  				 Inst3&30180&29850&990&360\\  
\hline    
  				 Inst4&26230&25780&9880&1010\\  
\hline    
  				 Inst5&28340&27200&7430&780\\  
\hline 
  			\end{tabular}
  			
  		}
  		
  	\end{center}                                                
  	
  \end{table*}



\newpage
\begin{thebibliography}{00}
	
\bibitem{ simpl 1}
Hosseini E, Reinhardt L, Rawat DB. Optimizing Gradient Methods for IoT Applications. IEEE Internet of Things Journal. 2022 Jan 11.

\bibitem{ simpl 2}
Hosseini E, Reinhardt L, Ghafoor KZ, Rawat DB. Implementation and Comparison of Four Algorithms on Transportation Problem. InInternational Summit Smart City 360° 2022 (pp. 422-433). Springer, Cham.

\bibitem{ simpl 3}
Hosseini E, Ghafoor KZ, Emrouznejad A, Sadiq AS, Rawat DB. Novel metaheuristic based on multiverse theory for optimization problems in emerging systems. Applied Intelligence. 2021 Jun;51(6):3275-92.

\bibitem{ simpl 4}
Hosseini E. Cost-Flow Summation Algorithm Based on Table Form to Solve Minimum Cost-Flow Problem. arXiv preprint arXiv:2101.01103. 2020 Dec 29.

\bibitem{ simpl 5}
Hosseini E, Ghafoor KZ, Sadiq AS, Guizani M, Emrouznejad A. Covid-19 optimizer algorithm, modeling and controlling of coronavirus distribution process. IEEE Journal of Biomedical and Health Informatics. 2020 Jul 28;24(10):2765-75.

\bibitem{ simpl 6}
Hosseini E, Sadiq AS, Ghafoor KZ, Rawat DB, Saif M, Yang X. Volcano eruption algorithm for solving optimization problems. Neural Computing and Applications. 2021 Apr;33(7):2321-37.

\bibitem{ simpl 7}
Ghafoor KZ, Kong L, Rawat DB, Hosseini E, Sadiq AS. Quality of service aware routing protocol in software-defined internet of vehicles. IEEE Internet of Things Journal. 2018 Oct 11;6(2):2817-28.

\bibitem{ simpl 8}
Hosseini E. Presentation and solving non-linear quad-level programming problem utilizing a heuristic approach based on Taylor theorem. Journal of Optimization in Industrial Engineering. 2018 Mar 1;11(1):91-101.

\bibitem{ simpl 9}
Hosseini E. Three new methods to find initial basic feasible solution of transportation problems. Applied Mathematical Sciences. 2017;11(37):1803-14.

\bibitem{ simpl 10}
Hosseini E. Solving linear tri-level programming problem using heuristic method based on bi-section algorithm. Asian J. Sci. Res. 2017;10:227-35.

\bibitem{ simpl 11}
Hosseini E. Laying chicken algorithm: a new meta-heuristic approach to solve continuous programming problems. J Appl Comput Math. 2017;6(1):1-8.

\bibitem{ simpl 12}
Hosseini E. Big bang algorithm: A new meta-heuristic approach for solving optimization problems. Asian Journal of Applied Sciences. 2017;10(3):134-44.

\bibitem{ simpl 13}
Hosseini E, Kamalabadi IN, Daneshfar F. Solving Non-Linear Bi-Level Programming Problem Using Taylor Algorithm. InHandbook of Research on Modern Optimization Algorithms and Applications in Engineering and Economics 2016 (pp. 797-810). IGI Global.

\bibitem{ simpl 14}
Hosseini E, Kamalabadi IN. Line search and genetic approaches for solving linear tri-level programming problem. Int J Manag Acc Econ. 2015;1(4).

\bibitem{ simpl 15}
Hosseini E, Nakhai Kamalabadi I. Two approaches for solving non-linear bi-level programming problem. Advances in Research. 2015;3(5):512-25.

\bibitem{ simpl 16}
Hosseini E, Kamalabadi IN. Bi-section algorithm for solving linear Bi-level programming problem. Int. J. Sci. Eng. 2015;1:101-7.

\bibitem{ simpl 17}
Hosseini E, Kamalabadi IN. Smoothing and solving linear quad-level programming problem using mathematical theorems. Int. J. Math. Comput. Sci. 2015;1:116-26.

\bibitem{ simpl 18}
Hosseini E, Kamalabadi IN. A modified simplex method for solving linear-quadratic and linear fractional bi-level programming problem. Global Journal of Advanced Research. 2015;1(2):142-54.

\bibitem{ simpl 19}
Nakhai Kamalabadi I, Hosseini E, Fathi M. Enhancing the solution method of linear Bi–level programming problem based on enumeration method and dual method. Journal of Advanced Mathematical Modeling. 2014 Aug 23;4(1):27-53.

\bibitem{ simpl 20}
Hosseini E, Kamalabadi IN. Line search and genetic approaches for solving linear tri-level programming problem. Int J Manag Acc Econ. 2015;1(4).

\bibitem{ simpl 21}
Hosseini E, Kamalabadi IN. Solving linear bi-level programming problem using two new approaches based on line search and taylor methods. Int J Manage Sci Education. 2014 Nov;2(6):243-52.

\bibitem{ simpl 22}
Hussein E, Kamalabadi I. Taylor approach for solving nonlinear bilevel programming problem. Adv. Comput. Sci., Int. J. 2014;3(5):91-7.

\bibitem{ simpl 23}
Hosseini E, Kamalabadi IN. Solving Linear-Quadratic Bi-Level Programming and Linear-Fractional Bi-Level Programming Problems Using Genetic Based Algorithm. Applied Mathematics and Computational Intelligence. 2013;2.	


	
\end {thebibliography}
\end{document}